\documentclass[12pt]{amsart}
\usepackage[latin1]{inputenc}
\usepackage{amsmath}
\usepackage{amsfonts}
\usepackage{amssymb}
\usepackage{amsmath,bbm,amssymb,amsxtra}
\usepackage{enumerate}
\allowdisplaybreaks
\usepackage{epsfig}
\usepackage{ae}

    \newcommand{\refeq}[1]{(\ref{#1})}

\newtheorem{lemma}{Lemma}[section]
\newtheorem{proposition}[lemma]{Proposition}

\newtheorem{theorem}[lemma]{Theorem}
\newtheorem{corollary}[lemma]{Corollary}
\newtheorem{remark}[lemma]{Remark}
\newtheorem{definition}[lemma]{Definition}

\newcommand{\prop}[1]{\begin{proposition}\label{#1}
\sl }
\newcommand{\eprop}{\end{proposition}}
\newcommand{\thm}[1]{\begin{theorem}\label{#1}
\sl }
\newcommand{\ethm}{\end{theorem}}
\newcommand{\nn}{\nonumber}

\newcommand{\lem}[1]{\begin{lemma}\label{#1}
\sl }
\newcommand{\elem}{\end{lemma}}

\newcommand{\defin}[1]{\begin{definition}\label{#1}
\sl }
\newcommand{\edefin}{\end{definition}}
\newcommand{\eps}{\epsilon}
\newcommand{\beqno}{\begin{eqnarray*}}
\newcommand{\eeqno}{\end{eqnarray*}}
\newcommand{\beqla}[1] {\begin {eqnarray}\label{#1}}
\def\eeq {\end {eqnarray}}
\newcommand{\beq}{\begin {eqnarray}}
\newcommand{\lip}{\operatorname{Lip}}

\newcommand{\Wk}{\operatorname{Wk}}
\newcommand{\diam}{\operatorname{diam}}
\newcommand{\R}{{\mathbb R}}
\newcommand{\midd}{{\; :\;}}

\newcommand{\Z}{{\mathbb Z}}
\newcommand{\N}{{\mathbb N}}

\newcommand{\ra}{\rightarrow}
\newcommand{\sub}{\subset}
\renewcommand{\:}{\colon}

\newcommand{\no}{\noindent}
\newcommand{\trace}{{\mathbb T\hskip-4pt{\mathbb R}} }

\author[Bonk]{Mario Bonk}
\address{Department of Mathematics, University of California, Los Angeles,  Box 95155, Los Angeles, CA, 90095-1555, USA}
\email{mbonk@math.ucla.edu}
\title[Sobolev spaces and hyperbolic fillings]{Sobolev spaces and hyperbolic fillings} 
\author[Saksman]{Eero Saksman}
\address{Department of Mathematics and Statistics, University of
Helsinki, PO~Box~68, FI-00014 Helsinki, Finland}
\email{eero.saksman@helsinki.fi}
\date{March 28, 2015}
\thanks{M.B.\ was supported by NSF grants DMS-0456940, DMS-0652915, DMS-1058283, DMS-1058772, and DMS-1162471}
\thanks{E.S.\  was
supported by the Finnish CoE in Analysis and Dynamics Research,
and by the Academy of Finland, projects 
113826  and  118765}
\begin{document}
\maketitle
\begin{abstract} Let $Z$ be  an Ahlfors $Q$-regular  compact metric measure space, where $Q>0$. For $p>1$ we introduce a new (fractional)  Sobolev space $A^p(Z)$ consisting of functions whose
extensions to the hyperbolic filling of $Z$ satisfy    a weak-type 
gradient condition. If $Z$ supports     a $Q$-Poincar\'e inequality with $Q>1$,  then  
$A^{Q}(Z)$ coincides with the familiar (homogeneous) Haj\l asz-Sobolev space.  \end{abstract}

\section{introduction}\label{se:intro}
In this paper we consider Ahlfors $Q$-regular compact metric measure spaces $Z= (Z,d , \mu)$, where $Q>0$.  We are interested in a certain Sobolev space on $Z$ and 
its   relation  to a function space that can be defined on a suitable hyperbolic filling of $Z$. Our results complement earlier work by  Bourdon and Pajot \cite{BP}, and by Connes, Semmes, Sullivan and
Teleman \cite[Appendix]{CST}.
In order to state our main theorems,  we first  have to discuss some basic concepts and set up some  notation. More details can be found in the later sections.

The hyperbolic filling $X$ of $Z$ (see Section~\ref{s:Fill}) is a  simplicial  graph $X=(V,E)$ with  vertex set $V$  and edge set $E$. It carries a natural path metric obtained by identifying each edge $e$ in $X$ with a copy of the unit interval.  Equipped with this  metric,  $X$  is Gromov hyperbolic and its  boundary at infinity $\partial_\infty X$  can be identified with $Z$. The construction of  $X$ depends on some  choices,  but  $X$ is uniquely determined up to quasi-isometry.  

In our construction,  the vertex set $V$ is given by a collection of metric balls $B$ in $Z$. Then   for an integrable function $f\in L^1(Z)$  (with $\mu$ being the underlying measure on $Z$) one can define the {\em Poisson extension} $u=Pf\: V\ra \R$ by setting 
$$
u(B)=\frac {1}{\mu(B)}\int_B f\, d\mu 
$$
for $B\in V$. 
 For each edge $e\in E$ we  choose one of the two  vertices incident with $e$ as the initial point $e_-$  and the other vertex as the terminal point $e_+$ of $e$.  
If $u\: V \ra \R$ is a function on $V$, we can then define the {\em gradient} $du\: E\ra \R$  of $u$ 
as   
\beqla{gradient} du(e)=u(e_+)-u(e_-)
\eeq
for $e\in E$. Note that both operators $f\mapsto Pf$ and $u\mapsto du$ are linear. 

If $M$ is a countable set and $p\ge 1$, then we denote by $\ell^{p,\infty}(M)$  the {\em weak-type $\ell^p$-space} consisting of all functions $s\: M \ra \R$ such that there exists a constant $C\ge 0$ with
$$ \# \{ m\in M: |s(m)|>\lambda\} \le\biggr(\frac{C}{\lambda}\biggl)^p$$
for all $\lambda>0$. For $p>1$ one can find  a norm $\Vert s\Vert_{\ell^{p,\infty}}$ for such   functions $s$ that is comparable  to  the infimum over all constants $C$ in the  previous inequality (see  Section~\ref{s:seq}).

Given these definitions one can define a  Sobolev space  as follows. 
\begin{definition}\label{de:QH} Let  $p>1$. Then  the Sobolev space $\dot A^p(Z)$ consists  of all  functions $f\in L^1(Z)$ such that the Poisson extension $u=Pf$ satisfies $du\in\ell^{p,\infty}(E)$.  
A semi-norm on this space is obtained  by setting 
\beqla{Anorm}
\| f\|_{\dot A^p}:=\|d(Pf)\|_{\ell^{p,\infty}} .
\eeq
\end{definition}
The semi-norm in \eqref{Anorm} does not distinguish functions that differ by an additive constant; so $\dot A^p(Z)$ is a Sobolev space of homogeneous type. 
 To promote the expression in  \eqref{Anorm} to a norm, we set
$A^p(Z)=\dot A^p(Z)/\R$, where $\R$ stands for the space of (almost everywhere)  constant functions 
on  $Z$. Strictly speaking,  the elements in $A^p(Z)$ are equivalence classes  in $L^1(Z)$ modulo constant functions, but we prefer to represent an  element in $A^p(Z)$ as an integrable function $f$  in $L^1(Z)$, defined up to an additive constant and up to changing the function on a set of measure zero. We  set 
$\| f\|_{A^p}:=\|d(Pf)\|_{\ell^{p,\infty}},$ which is well-defined on $A^p(Z)$.

\begin{proposition}\label{co:subspace} Let $p>1$. Then the map $f\in A^p(Z)\mapsto \| f\|_{A^p}$ defines a norm on $A^p(Z)$. With this norm, 
 $A^p(Z)$ is a Banach space that is isomorphic to
a closed subspace of $\ell^{p,\infty}(E).$
\end{proposition}

By definition the space $A^p(Z)$ depends a priori  on the choice of the hyperbolic filling $X=(V,E)$ of $Z$, but  we will see that $A^p(Z)$ is independent of the choice of $X$
(see Corollary~\ref{cor:indfill} and Remark~\ref{rem:fill}).

One can  find an isomorphism of  the  Sobolev space $A^p(Z)$ with a  space that can be constructed from  functions  on the hyperbolic filling $X=(V,E)$. 
For this we fix $p> 1$ and  define an equivalence relation on real-valued functions
 on $V$ as follows: if $u,u'\: V \ra \R$, we write $u\sim u'$ if there exists a constant $c\in \R$ such that  $u-u'-c\in \ell^{p,\infty}(V)$.  It is clear that  $\sim$ is an equivalence relation and we denote by $[u]$  the corresponding equivalence class of a function $u\: V \ra \R$. We now define $\Wk^{1,p}(X)$ as the real vector space  of all equivalence classes $[u]$ 
 of functions $u\: V\ra \R$ with $du\in \ell^{p,\infty}(E)$. So 
\beqla{wkl^qco}
\Wk^{1,p}(X)=\{ u\: V\ra \R:  du\in  \ell^{p,\infty}(E)\}/(\R+\ell^{p,\infty}(V)),
\eeq
 where $\R$ represents the space of constant functions on $V$. Note that 
 if $u\in \R+\ell^{p,\infty}(V)$, then  $du\in  \ell^{p,\infty}(E)$ (see Section~\ref{s:Fill}). 
 The space $\Wk^{1,p}(X)$ carries a natural semi-norm given by 
  \beqla{semnordefWk} \Vert u\Vert_{\Wk^{1,p}}:= \inf\{ \Vert du'\Vert_{\ell^{p,\infty}}: u'\in [u]\}. 
\eeq
Note that by definition $\Vert u\Vert_{\Wk^{1,p}}$ only depends on $[u]$ and not on the chosen representative $u$ in $[u]$. So $\Vert u\Vert_{\Wk^{1,p}}$ really defines a semi-norm on $\Wk^{1,p}(X)$, but for simplicity we suppress the distinction between $u$ and $[u]$ in our notation for the semi-norm. 

It is easy to see  that under some mild assumptions  the space  $\Wk^{1,p}(X)$ is  invariant under quasi-isometries of $X$ (see Proposition~\ref{prop:qiinv}).

 If  $p>1$, and  $u\: V\ra \R$ is  a function with $du\in \ell^{p,\infty}(E)$, then one can show that 
 it has a well-defined {\em trace} $\trace u\in L^1(Z)$ on $Z=\partial_\infty X$ 
 (see Lemma~\ref{le:trace}).

We can now state the main result of this paper.
\begin{theorem}\label{main2}
Let $p>1$, and $Z$ be an Ahlfors $Q$-regular space, where $Q>0$. Then the expression   $\refeq{semnordefWk}$ defines a norm on   $\Wk^{1,p}(X)$. Equipped with this norm,  $\Wk^{1,p}(X)$ is a Banach space. 

Moreover, 
 the map 
 $[u]\mapsto \trace u$  is an isomorphism between the Banach spaces $\Wk^{1,p}(X)$ and $A^p(Z)$ with the inverse obtained from  the Poisson extension by $f\mapsto [Pf]$.  
\end{theorem} 

For general Ahlfors regular spaces the relation of our space $A^p(Z)$ to other known Sobolev-type spaces it not clear. It is easy to see that the (homogeneous)
 fractional Haj\l asz-Sobolev $\dot M^{\alpha,p}(Z)$ is always contained in $\dot A^p(Z)$ for $\alpha=Q/p$ (see Section~\ref{s:PI}).
 If  $Z$ satisfies a Poincar\'e inequality, then 
 our space can be identified  with the   Haj\l asz-Sobolev space for the endpoint case $p=Q$.

\begin{theorem}\label{main} Let $Z$ be  an Ahlfors $Q$-regular compact metric measure space   that supports a  $Q$-Poincar\'e inequality, where $Q>1$. Then we have $\dot A^{Q}(Z)=\dot M^{1,Q}(Z)$, with comparability of semi-norms.
\end{theorem}

One can also show that under the assumptions of this theorem the space 
$\dot A^p(Z)$ consists only of (almost everywhere) constant functions for $1< p<Q$ (see Proposition~\ref{prop:consfunc}). 

A  reader with an orientation towards more classical analysis may  wonder about analogs of  our results in a standard Euclidean setting.   
We have included a discussion on this in Section~\ref{s:rn}. Here we will relate 
integrability  properties of a  function $f$ in $\R^n$ with  
weak-type conditions for the classical Poisson extension $u$ of $f$ defined on upper half-space 
$\R_+^{n+1}$ (see Proposition~\ref{prop:fuweak}). This can be used to characterize functions $f$ 
with a distributional gradient in $L^n(\R^n)$ in terms of a weak-type condition of the hyperbolic gradient $\nabla_hu$;  see Corollary~\ref{co:Q-sobolevweak} which corresponds to Theorem~\ref{main}.  

 The definition \eqref{wkl^qco} of our space $\Wk^{1,p}(X)$ is similar to  a definition 
 in \cite{BP}, where the requirement is in terms of an $\ell^p$-condition instead of a weak-type condition. The   space obtained in this way admits  an identification with the  $\ell^p$-cohomology 
 of the hyperbolic filling $X$ in degree $1$ and can  also  be identified with a Besov space $B_p(Z)$ on the boundary $Z=\partial_\infty X$. In \cite{BP} it was also shown 
 that in the setting as in Theorem~\ref{main}, the space $B_p(Z)$ is  interesting only for $p>Q$, because it is trivial and  consists  of constant functions  in the endpoint case $p=Q$. Theorem~\ref{main} suggests that one should modify 
 the definition of $\ell^p$-cohomology if one wants to obtain interesting function spaces for critical exponents. Namely, one can set up a cohomology theory for (infinite) simplicial complexes (satisfying additional natural geometric conditions), where the requirement is that cochains belong to $\ell^{p,\infty}$ instead of $\ell^p$. Our notation for the space $\Wk^{1,p}(X)$ is suggested by the fact that it represents  weak-type $\ell^p$-cohomology in  degree $1$.  
 Many important features of $\ell^p$-cohomology such as quasi-isometric invariance properties remain  valid in this context of  weak-type $\ell^p$-cohomology. 
 This is an interesting direction to pursue, but we will not do this in this paper. 
 
 Our paper is organized as follows. In Section~\ref{s:seq} we review some facts 
 about sequence spaces and weak-type conditions. We consider the space $\Wk^{1,p}(X)$ on general simplicial graphs $X$ in Section~\ref{s:Fill}, where we prove a quasi-isometric invariance property. In  Section~\ref{s:Ahl} we specialize to Ahlfors regular spaces $Z$ and their hyperbolic fillings. There  we will prove Proposition~\ref{co:subspace} and Theorem~\ref{main2}. In Section~\ref{s:PI} we discuss  Sobolev spaces and their relation to our space $\dot A^p(Z)$. In particular, we  give a proof of Theorem~\ref{main}. As was already mentioned, Section~\ref{s:rn} is devoted to results about the classical Poisson extension in $\R^n$.

\medskip 
\no
{\bf Acknowledgement.} The authors are indebted to Marc Bourdon, Bruce Kleiner, Pierre Pansu, and Tomas Soto  for many interesting  discussions relating to the topic of this paper. They also thank Jeff Lindquist for a careful reading of a draft of the paper.

\section{Sequence spaces}\label{s:seq} 
In the following, we will extensively use the notation  $A\lesssim B$, $A\gtrsim B$, and $A\simeq B$ for quantities $A$ and $B$ to indicate the existence of an
implicit  constant $C\ge 1$ depending on some inessential  parameters such 
$A\le CB$, $A\ge B/C$ and $A/C\le B\le CA$, respectively.  Of course,  it depends  on the context
which parameters can be safely  ignored. We will mostly leave it to the judicious readers to make their own judgements about this. 

For  $p\in [1,\infty)$  we denote 
by $\ell^p$ the space of all real-valued sequences $s=\{x_n\}\in \R^\N $ such that 
$$\Vert s\Vert_{\ell^p}:=\biggl( \sum_{n=1}^\infty |x_n|^p\biggr)^{1/p}<\infty. $$
Note that $\Vert s\Vert_{\ell^q}\le \Vert s\Vert_{\ell^p}$  and so 
$\ell^p\sub \ell^q$ for $1\le p\le q$.  

We denote by $\ell^{p,\infty}$ the space of all  sequences $s=\{x_n\} \in \R^\N $
for which
there exists a constant $C\ge 0$ with 
 \beqla{weaknonnorm} 
\# \{ n\in \N: |x_n|>\lambda\} \le \biggl(\frac{C}{\lambda}\biggr)^p
\eeq for all $\lambda>0$. Here $\#M\in \N_0\cup\{\infty\}$ denotes the cardinality of a set $M$. 
 Note that $\ell^p\subset \ell^{p,\infty}$.  

For a given sequence $s=\{x_n \}\in \R^\N$ we denote the infimum of all constants $C$ for which \eqref{weaknonnorm} is valid as 
$ \Vert s\Vert^*_{\ell^{p,\infty}}.$  
We then have the homogeneity property  $ \Vert as\Vert^*_{\ell^{p,\infty}}=|a| \cdot  \Vert s\Vert^*_{\ell^{p,\infty}}$ for $a\in \R$, but we  do not have subadditivity. So $ \Vert s\Vert^*_{\ell^{p,\infty}}$  does not define  a norm on $\ell^{p,\infty}$. For $p>1$ we consider the expression 
\beqla{weaknorm} 
\| s\|_{\ell^{p,\infty}}:=\sup\big\{ n^{-1+1/p}\sum_{j\in E}|x_j|: n\in \N,\, E\subset \N,\,  \#E=n\big\}.
\eeq
The quantities
 $ \Vert s\Vert^*_{\ell^{p,\infty}}$ and $ \Vert s\Vert_{\ell^{p,\infty}}$ are comparable  (see Lemma~\ref{le:weaktostrong}~(i) below)  and it is easy to verify  that $\| s\|_{\ell^{p,\infty}}$ defines a norm on $\ell^{p,\infty}$. We will equip $\ell^{p,\infty}$ with this norm in the following, but we will freely switch between $\| s\|_{\ell^{p,\infty}}$ and the comparable expression $ \Vert s\Vert^*_{\ell^{p,\infty}}$ whenever convenient. 
The space   $\ell^{p,\infty}$ is  a Banach space, and one can  show that it is non-reflexive (see \cite[Section~1.4]{Gr} for general background on  
weak $L^p$-spaces).

We will also consider sequences 
 indexed by more general countable index sets $M$, 
 or, more precisely, functions $s\: M\ra \R$. For clarity we will then denote the corresponding spaces  by $\ell^p(M)$ and   $\ell^{p,\infty}(M)$, but will suppress $M$ in the notation for norms. 


We need  several simple observations. 
\begin{lemma}\label{le:weaktostrong} 
\begin{itemize}

\item[\textnormal{(i)}] For each $p>1$ there  exists a constant $C\ge 1$ such that 
$$ \|s\|^*_{\ell^{p,\infty}}\le \|s\|_{\ell^{p,\infty}}\le C
\|s\|^*_{\ell^{p,\infty}}
$$ 
for all $s\in \R^\N$. 

\smallskip
\item[\textnormal{(ii)}]
 Let $p>1$ and $r\in [1,p)$. Then  there exists a constant $C'>0$ such for every sequence   
 $s=\{x_1,\ldots ,x_m,0,0,\ldots\}\in \R^{\N}$ with $m\in \N$ we have 
 $$
\| s\|_{\ell^p}\leq C'(1+\log m)^{1/p}\|s\|_{\ell^{p,\infty}}\quad {\rm and}\quad \| s\|_{\ell^r}\leq C'm^{1/r-1/p} \|s\|_{\ell^{p,\infty}}.
$$
\end{itemize} 
\end{lemma}
\begin{proof} (i) Let $s=\{x_n\}\in \R^\N$ be arbitrary.

To show the first inequality, we may assume 
$\|s\|_{\ell^{p,\infty}}=1$. If $\lambda>0$ and $E\sub \N$ is a finite set 
such that $|x_j|>\lambda$ for $j\in E$,  we then  have 
$$\lambda(\#E) \le \sum_{j\in E} |x_j|\le (\#E)^{1-1/p}. $$ 
It follows that $\#E\le 1/\lambda^p$.  This implies $\|s\|_{\ell^{p,\infty}}^*\le 1$
as desired. 

To establish  the second inequality, we may assume that $\|s\|_{\ell^{p,\infty}}^*=1$. 
Let $E\sub \N$ be an arbitrary finite non-empty set and set $m=\#E$. 
Then for each $k\in \Z$ the  number of elements $j\in E$ 
  with $|x_j| \in (2^{-k},2^{1-k}]$ is bounded above by $\min \{m, 2^{pk}\}$.
Hence 
$$ \sum_{j\in E} |x_j| \le  \sum_{k=0}^\infty 2^{1-k}  \min \{m, 2^{pk}\}\le Cm^{1-1/p}.
$$ 
Here $C\ge 1$ can be chosen only  depending  on $p$. It follows that 
$\|s\|_{\ell^{p,\infty}}\le C$ as desired.  

\smallskip 
(ii) Let $s\in \R^\N$ be as in the statement.  
 We may assume $\|s\|_{\ell^{p,\infty}}=1$ which implies $\|s\|_{\ell^{p,\infty}}^*\le 1$ by (i). Then, similarly as in the proof of (i),   for each $k\in \Z$ the  number of elements $x_j$ in the sequence $s$  with $|x_j| \in (2^{-k},2^{1-k}]$ is bounded above by $\min \{m, 2^{pk}\}$. Hence 
$$ \| s\|_{\ell^p}^p=\sum_{j=1}^m |x_j|^p\lesssim \sum_{k=0}^\infty 2^{-pk}  \min \{m, 2^{pk}\}\lesssim 1+\log m.
$$
The other estimate follows from a similar computation.
\end{proof}

\begin{lemma}\label{le:weakbdd} Let $p> 1$ and 
$E\subset \N \times \N$. For $k,n\in \N$ define 
 $$E_n=\{k\in \N: (n,k)\in E\} \text{ and }
E^k= \{n\in \N: (n,k)\in E\},$$  and suppose that for some $N\in \N$ 
we have  $\#E_n\le N$ and $ \#E^k\le N$ for all $k,n\in \N$. Then there exists a constant $C=C(p,N)>0$ with the following property:

If  $s=\{x_n\}$ and $t=\{y_n\}$ are sequences in  $\R^\N$ such that 
$$ |x_n|\le \sum_{k\in E_n}|y_k|$$
for all $n\in \N$,  
then $ \Vert s\Vert_{\ell^{p,\infty}}\le C(p,N) \Vert t\Vert_{\ell^{p,\infty}}$.  
 \end{lemma}
\begin{proof} We may assume $\|t\|_{\ell^{p,\infty}}^*=1.$  Let $\lambda>0$
and suppose that $|x_n|>\lambda$. Then there exists $k\in E_n$ such that 
$|y_k|>\lambda/N$. Each such $k$ lies in at most $\#E^k\le N$   sets $E_n$.
Hence 
$$\#\{ n\in \N: |x_n|>\lambda\}\le N \#\{ k\in \N: |y_k|>\lambda/N\}\le N^{p+1}/\lambda^p
. $$
 The claim follows. 
\end{proof}

In  the next lemma it is convenient to view a sequence $s\in \R^\N$ as a function $s\: \N \ra\R$. 

\begin{lemma}\label{le:weakconv} Let $p\ge 1$ and suppose  that $s$ and $s_k$ for $k\in \N$ are sequences in $\R^\N$ such that for each $n\in \N$ we have $s_k(n)\to s(n)$ as $k\to \infty$. Then
$$ \Vert s\Vert^*_{\ell^{p,\infty}}\le \liminf_{k\to \infty}  \Vert s_k\Vert^*_{\ell^{p,\infty}}. $$ 

In particular, if  $\Vert s_k\Vert_{\ell^{p,\infty}}$ is uniformly bounded for $k\in \N$, then 
$s\in \ell^{p,\infty}$. 
 \end{lemma} 
\begin{proof} If $C_0:=\liminf_{k\to \infty}  \Vert s_k\Vert^*_{\ell^{p,\infty}}=\infty$, there is nothing to prove. So suppose  that $C_0<\infty$. Pick $C>C_0$. By passing to a subsequence if necessary, we may assume that 
 $\Vert s_k\Vert^*_{\ell^{p,\infty}}<C$ for all $k\in \N$. Let $\lambda>0$ be arbitrary. If 
 $|s(n)|>\lambda$ for some $n\in \N$, then $|s_k(n)|>\lambda$ for all sufficiently large $k$. Hence 
 $$ \#\{ n\in \N : |s(n)|>\lambda\} \le \liminf_{k\to \infty} \#\{ n\in \N : |s_k(n)|>\lambda\}
 \le (C/\lambda)^p,$$ and so  $\Vert s\Vert^*_{\ell^{p,\infty}}\le C$. If we let $C\to C_0$, the claim follows. 
\end{proof}

\section{Function spaces on simplicial graphs}\label{s:Fill} 
Before we turn to hyperbolic fillings of Ahlfors regular spaces, we will discuss  some facts for 
general simplicial graphs
 $X=(V,E)$.  We assume 
 that $X$ is connected and carries a path   metric obtained by identifying each edge $e\in E$ with a copy of the unit interval $[0,1]$. We will also  assume that each pair of vertices in $V$ is joined by at most one edge in $E$ and that  the degree of each vertex $v\in V$, i.e., the number of edges incident with $v$,  is uniformly  bounded from above.  This implies that if $R>0$,  then the number of edges and vertices contained in a ball $B\subset X$ of radius $R$ is uniformly bounded above only depending on $R$.
 
  As already mentioned in the introduction, it is convenient  choose one of the vertices incident with $e\in E$ as the initial point $e_-$  and the other  as the terminal point $e_+$ of $e$. For a function   $u\: V\ra \R$, we define the gradient $du\: E\ra \R$ of $u$ as in \eqref{gradient}.
 Then for  each $e\in E$ we have 
 $$|du(e)|=|u(e_+)-u(e_-)|\le |u(e_+)|+|u(e_-)|. $$
 Here  each vertex $v\in V$ appears as an endpoint $e_+$ or $e_-$  only for a uniformly bounded number of edges $e$. Hence  Lemma~\ref{le:weakbdd} implies that if $p> 1$, then 
 \beqla{lqdlq}
 \Vert du\Vert_{\ell^{p,\infty} } \le C \Vert u\Vert_{\ell^{p,\infty} } 
 \eeq 
with a constant $C>0$ independent of $u$. In other words,  the map $u\mapsto du$ is a bounded linear operator $d\: \ell^{p,\infty}(V)\ra  \ell^{p,\infty}(E)$. 
In particular, $\R+ \ell^{p,\infty}(V)$ is is a subspace of 
$\{ u\: V\ra \R: du\in \ell^{p,\infty}(E)\}$ and one can 
define $\Wk^{1,p}(X)$ as in \eqref{wkl^qco} for each $p> 1$. 
The space $\Wk^{1,p}(X)$ carries the semi-norm as defined in \eqref{semnordefWk}. 

Two metric spaces $(X,d_X)$ and $(Y,d_Y)$ are called {\em quasi-isometric} if there exists a map $\varphi\: X\ra Y$ such that for some constants $\lambda\ge1$ and $K\ge 0$ we have $$ \frac1\lambda d_X(x,x')-K\le d_Y(\varphi(x), \varphi(x'))\le \lambda d_X(x,x')+K$$ 
for all $x,x'\in X$, and 
$$ \sup_{y\in Y} \inf_{x\in X} d_Y(y, \varphi(x))\le K. $$
A  map $\varphi\: X\ra Y$  as in this definition is called a 
{\em quasi-isometry}.

 \begin{proposition} \label{prop:qiinv}
 Suppose $X=(V,E)$ and $X'=(V',E')$ are simplicial metric graphs that are connected and have uniformly bounded vertex degree. If $X$ and $X'$ are quasi-isometric, then 
  for each $p>1$ the spaces $\Wk^{1,p}(X)$ and $\Wk^{1,p}(X')$ are isomorphic.  
  \end{proposition}
  
\begin{proof}  We will show that there exists a linear bijection between the function spaces that gives comparability of the semi-norms.

Our assumptions imply that  there exist quasi-isometries   $\psi\: V'\ra V$ and 
$\varphi\: V\ra V'$ that are  coarse inverses of each other; more precisely, $\psi\circ \varphi$ and $\varphi\circ \psi$ are in bounded distance to the identity map on $V$ and $V'$, respectively. 

 Then any two preimages of a vertex $v\in V$ under $\psi$  have uniformly bounded distance. Since the vertex degree of $X$ is uniformly bounded, we conclude  that there exists a constant $N\in \N$ such that 
  $\#\psi^{-1}(v)\le N$ for all   $v\in V$. This implies that   if $u\in \R+\ell^{p,\infty}(V)$, then $u\circ \psi \in \R+\ell^{p,\infty}(V')$. 
  
  Suppose $u\: V \ra \R$ and $du\in \ell^\infty(E)$. If $e'\in E'$, then the vertices 
  $e'_+$ and  $e'_-$  are adjacent in $V'$. Then they have distance $1$ in $X'$, and so  their  
  images $\psi(e'_+)$ and $\psi(e'_-)$ have uniformly bounded distance in  $V$. We join these images by a geodesic segment, and estimate their function value difference by the triangle inequality. It then follows that  for  a radius  $R_0>0$ independent of $e'$ we have
  \begin{equation}\label{sumwk} |d(u\circ \psi)(e')|= |(u (\psi(e'_+))-u(\psi(e'_-))|\le 
  \sum_{e\subset  B(\psi(e'_+), R_0)} |du(e)|,
  \end{equation}
  where $B(a,r)\sub X$  denotes the open ball of radius $r>0$ centered at $a\in X$. 
  
 The number of the edges $e$ contributing  to  the last sum is uniformly bounded independent of $e'$; moreover, if a given edge   $e\in E$ contributes, then $e_+$ and $ \psi(e'_+)$ have uniformly 
 bounded distance. If we apply the coarse inverse $\varphi$ of $\psi$ here, we see that 
 $\varphi(e_+)$ and $ e'_+$ have uniformly bounded distance; so for given $e\in E$ there is  uniformly bounded  number of edges $e'\in E'$ so that $e$ contributes to the sum
 \eqref{sumwk}. These considerations show that we can apply Lemma~\ref{le:weakbdd} and we conclude that 
  $$ \Vert d(u\circ \psi)\Vert_{\ell^{p,\infty}}\le C \Vert du\Vert_{\ell^{p,\infty}}.
  $$ with a constant $C\ge 0$ independent of $u$.

 By what we have seen,  the linear map 
 $$[u]\in \Wk^{1,p}(X) \mapsto [u\circ \psi]\in \Wk^{1,p}(X')$$ is well-defined, and we get a uniform  semi-norm bound
  $$ \Vert  u\circ \psi\Vert_{\Wk^{1,p}} \le C \Vert u\Vert_{\Wk^{1,p}} . $$ 
  
Of course, the roles of $\varphi$ and $\psi$ are completely symmetric, and we can apply the previous considerations to $\varphi$. So  we get a well-defined linear map  $$[u']\in \Wk^{1,p}(X') \mapsto [u'\circ \varphi]\in \Wk^{1,p}(X)$$ with a corresponding semi-norm bound.
 To finish the proof, it is enough to show that  the maps $[u]\mapsto [u\circ \psi]$ and $[u']\mapsto [u'\circ \psi]$ are inverse to each other, or equivalently that $[u]=[u\circ\psi\circ \varphi]$ and $[u']=[u' \circ \varphi\circ \psi]$. By symmetry, it is enough to show that 
 $[u]=[u\circ\psi\circ \varphi]$ for $[u]\in  \Wk^{1,p}(X)$. 

So let  $u\: V \ra \R$ with $du\in \ell^{p,\infty}(E)$ be arbitrary. If $v\in V$, 
then $(\psi\circ \varphi)(v)\in B(v,R_1)$, where $R_1$ is independent of $v$.
Hence 
$$|u(v)-(u\circ\psi\circ \varphi)(v)|\le \sum_{e\subset B(v,R_1)} |du(e)|.$$
Again there is only a uniformly bounded number of edges $e$ contributing to this sum, and a given edge $e$ can only appear  for a uniformly bounded number of vertices $v\in V$.  
So by Lemma~\ref{le:weakbdd} we have 
$$ \Vert u-u\circ\psi\circ \varphi\Vert_{\ell^{p,\infty}}\lesssim \Vert du\Vert_{\ell^{p,\infty}}<\infty.
  $$
  In particular, $u-u\circ\psi\circ \varphi \in \ell^{p,\infty}(V)$, and so 
  $[u]=[u\circ\psi\circ \varphi]$ as desired. 
 \end{proof}

\section{Ahlfors regular spaces and hyperbolic fillings}\label{s:Ahl} 

In this and the following sections $Z= (Z,d , \mu)$ is  an Ahlfors $Q$-regular compact metric measure space, where $Q>0$. Here $d$ is a metric and $\mu$ a Borel  measure on $Z$. We will  assume that the diameter $\diam(Z)$ of $Z$ is equal to $1$.  This can always be achieved by a possible rescaling of the metric. 
We use both $\mu(E)$ and $|E|$ to denote the measure of a (Borel) set $E\subset Z$. 
 Ahlfors $Q$-regularity then  means that $|B|\simeq r^Q$ whenever  $B$ is a  ball of radius $r\le 1$ in $Z$.

 If $a\in Z$ and $r>0$, we denote by $B(a,r)$ the open ball in $Z$ of radius $r$ centered at $a$. If $\Lambda\ge 1$ and $B=B(a,r)$ is a ball, we set $\Lambda B=B(a,\Lambda r)$. 
For $p\ge 1$ we denote by $L^{p}(Z)$ the space of all measurable functions $f\: Z\ra \R$ such that 
$$ \Vert f\Vert_{L^p}:= \biggl( \int_Z |f|^p\, d\mu\biggr)^{1/p}<\infty. $$ 
If $f\in L^1(Z)$ and $B$ is a ball in $Z$, we define
$$ f_B=\frac{1}{|B|} \int_B f\, d\mu. $$

We will  now  review  the construction of the  hyperbolic filling of $Z$. 
We will mostly follow  \cite[Section 2.1]{BP} (see also  \cite[Chapter 6]{BS}). 
Note that we can apply the considerations  in \cite{BP}, because as an Ahlfors 
regular space,   $Z$ is {\em doubling} and {\em uniformly perfect}
(see \cite[Sections 10.13 and 11.1]{He} for this terminology).

For each $n\in \N$ we choose a maximal $2^{-n}$-separated set $Z_n$ of points in $Z$. Then   the balls $B(z, 2^{-n})$, $z\in Z_n$, form a cover 
of $Z$. The hyperbolic filling is a graph that essentially records the incidence relations of the balls in these covers. In order to get good properties it is necessary to  enlarge the radii of these balls by a factor 
$\Lambda>1$ (this was overlooked in 
\cite[Section 2.1]{BP}, and causes problems in the  proof of \cite[Lemme 2.2]{BP}, but the argument is valid if one uses  enlarged balls).  

Accordingly,  for $n\in \N$ we let $V_n$ 
 be the collection  of all balls $B=B(z, 2^{1-n})$ with  $z\in Z_n$, where we refer to $n$ as the {\em level} of the  ball $B$. Here we add 
 $B=Z$ as a  ball  on level $0$ and set $V_0=\{Z\}$. The vertex set $V$ of our graph is the (disjoint) union of the sets $V_n$, $n\in \N_0$.  We  join two distinct vertices as represented by balls $B$ and $B'$ 
by an edge  $e$ if their levels  differ by at most $1$ and if $B\cap B'\ne \emptyset$. We 
write $B'\sim B$ in this case and  denote by  $E$  the set of these edges $e$.  The simplicial graph $X=(V,E)$ carries a natural path metric obtained by  identifying each edge $e$ in $X$ with a unit interval.  Then the graph $X$ equipped with this  metric is Gromov hyperbolic and the boundary at infinity $\partial_\infty X$ of $X$ can be identified with $Z$. The graph $X$ is also connected and the degree of each vertex is uniformly bounded from above.  In particular, we can apply the considerations in Section~\ref{s:Fill} to  $X$. 

The construction of $X=(V,E)$ depends on choices, namely on  the sets $Z_n$. One can also choose  a scale sequence $\Lambda^{-n}$ with a different parameter $\Lambda>1$, instead of the sequence $2^{-n}$ (as in  \cite{BS} and \cite{BP}).  These  constructions lead to hyperbolic fillings  that are  quasi-isometric to our filling  $X$ of $Z$ (this can easily be proved directly or one can invoke a general fact such as  \cite[Th\'eor\`eme 2.3]{BP}). 
By Proposition~\ref{prop:qiinv} this quasi-isometric ambiguity is  irrelevant for the  function spaces $\Wk^{1,p}(X)$ we are interested in.  

%

For $n\in \N$ the set $V_n$ is the sphere of radius $n$ centered at $B_0=Z\in V$  in the subset $V$ of the hyperbolic filling $X$.  Since $Z$ is Ahlfors $Q$-regular, we have  
$\#  V_n\simeq 2^{Qn}.$

 We make the conventions underlying the definition in \eqref{gradient}. It is useful to also use  variants of the gradient $du$; namely, if $u\: V\ra \R$ is arbitrary, we define 
$  \tilde d u \: V\ra \R$ by  
 $$ \tilde d u(B)= \sum_{B'\in V,\, B'\sim B} |u(B')-u(B)| $$  
 for $B\in V$, and $  \tilde d_n u\: V_n\ra \R$ as the restriction of   $\tilde d u $ to $V_n$ for $n\in \N_0$. Note that in contrast to the gradient $du$, the function $ \tilde d u$ is defined on $V$ and not on $E$ and  $u\mapsto \tilde d u$ is not a linear map. 
 
 If $p\ge 1$ and $u\: V \ra \R$ is arbitrary, then for each $B\in V$ we have 
 $$   \tilde d u(B)\le \sum \{|du(e)|: \text{$e\in E$ and $B=e_+$ or $B=e_-$}\}. $$
By Lemma~\ref {le:weakbdd} this implies  $\Vert \tilde d u\Vert_{\ell^{p,\infty}}\lesssim  \Vert d u\Vert_{\ell^{p,\infty}}$. In the other direction, we 
have $|du(e)|\le \tilde du(e_-)$ for $e\in E$, which again by   Lemma~\ref {le:weakbdd}
implies that $ \Vert d u\Vert_{\ell^{p,\infty}} \lesssim   \Vert \tilde d u\Vert_{\ell^{p,\infty}}$.
It follows that 
\begin{equation}
\label{dduin}
 \Vert \tilde d u\Vert_{\ell^{p,\infty}}\simeq \Vert d u\Vert_{\ell^{p,\infty}} 
 \end{equation}
 with an implicit multiplicative constant independent of $u\: V\ra \R$.

For each $n\in \N$ the balls $\frac12 B$ with $B\in V_n$ form an open cover of 
$Z$. Accordingly,  we  can choose a Lipschitz partition of unity on $Z$ given by finitely many non-negative Lipschitz functions $\{\psi_B\}_{B\in V_n}$ on $Z$ 
such that  $$\sum_{B\in V_n}\psi_B=1, $$ 
 \begin{equation}\label{suppBincl}
  {\rm supp}(\psi_B):=\overline {\{ u\in Z: \psi(u)\ne 0\} }
 \subset B,
 \end{equation}
 and 
 \begin{equation}\label{Lipbdd} 
 \mathrm{Lip} (\psi_B):= \sup \biggl\{\frac { |\psi_B(u)-\psi_B(v)|}{d(u,v)}: u,v\in Z,\, u\ne v\biggr\}\lesssim2^n.
 \end{equation} 
See   \cite[p.~431, B.7.4 Lemma] {Se}) for the existence of such a partition of unity (the inclusion \eqref{suppBincl} which is stronger than in this reference can easily be obtained from a slight modification in the proof).
 
 For the ball $B=Z$ on level $0$ we choose $\psi_B=1$. 
 
For $u\: V\ra \R$ we define a function $T_nu\: Z \ra \R$ as 
 $$T_nu=\sum_{B\in V_n}u(B)\psi_B.
$$
Then  for $n\in \N_0$ we have 
\begin{align}\label{diffest}
|T_nu-T_{n+1}u|&\le \sum_{B\in V_n}\sum_{B'\in V_{n+1} }|u(B)-u(B')|\psi_B\psi_{B'}\\ 
&\le \sum_{B\in V_n}\tilde du(B)\psi_B. \nonumber
\end{align}
Note that the last inequality follows from the fact that $ \psi_B\psi_{B'}=0$ for $B\in V_n$ and $B'\in V_{n+1}$ unless $B\cap B'\ne \emptyset$ and so $B\sim B'$. 

Now   $\| \psi_B\|_{L^1}\lesssim 2^{-nQ}$ for $B\in V_n$, and so we obtain
\beqla{L1bdd}
\Vert T_nu-T_{n+1}u\Vert_{L^1}\lesssim  2^{-nQ} \Vert \tilde d_n u\Vert_{\ell^1}.  
\eeq
 The following lemma  provides an  inverse operation to the  Poisson extension.

\begin{lemma}\label{le:trace} Let $p>1$ and  $u\: V \ra \R$  with $du\in \ell^{p,\infty}(E).$ Then 
$$ \sum_{n=0}^\infty \|T_{n+1}u-T_nu\|_{L^1}\lesssim \| du\|_{\ell^{p,\infty}}.
$$ In particular, the  limit 
$$
\trace u:=\lim_{n\to \infty} T_nu.
$$ exists  both in $L^1(Z)$ and pointwise almost everywhere in $Z$. 
\end{lemma}

\begin{proof}  We have  $\#V_n\lesssim 2^{nQ}$,  and so  
 Lemma \ref{le:weaktostrong}~(ii) and \eqref{dduin} imply   that 
 $$\|\tilde d_nu\|_{\ell^1}\lesssim (2^{nQ})^{1-1/p}\| du\|_{\ell^{p,\infty}}.$$
Hence by \eqref{L1bdd}, 
\begin{align*}\sum_{n=0}^\infty \|T_{n+1}u-T_nu\|_{L^1}&\lesssim
\sum_{n=0}^\infty 2^{-nQ}
\|\tilde d_nu\|_{\ell^1} \\ &\lesssim \sum_{n=0}^\infty 2^{-nQ/p}\| du\|_{\ell^{p,\infty}}\lesssim 
\| du\|_{\ell^{p,\infty}}. 
\end{align*}   The claim easily  follows. 
\end{proof}

\begin{lemma}\label{lem:frecov} Let $p>1$.  If   $f\in L^1(Z)$ and  $d(Pf)\in \ell^{p,\infty}(E)$, 
 then $f=\trace (Pf)$. Moreover,
 $$ \biggl\Vert f-\frac{1}{|Z|}\int_Zf\,d\mu\biggr\Vert_{L^1}\lesssim \|d(Pf) \|_{\ell^{p,\infty}}. $$
\end{lemma}

\begin{proof}
To prove the first part, let $\epsilon>0$ be arbitrary. We know that $T_n(Pf)\to \trace (Pf)$ in $L^1(Z)$ as $n\to \infty$ by Lemma~\ref{le:trace}, and so $\Vert T_n(Pf)-\trace (Pf)\Vert_{L^1}<\eps$ for large $n$. 
 Since continuous functions are dense in $L^1(Z)$, we can find 
a continuous function $g$ on $Z$ with $\Vert g-f\Vert_{L^1}<\epsilon$.
It is clear that   $T_n(Pg)\ra g$ uniformly on $Z$  and hence also  
$\Vert T_n(Pg)-g\Vert_{L^1}<\eps$ for large $n$.

Note that if $h\in L^1(Z)$, then 
\begin{align*} \Vert T_n(Ph)\Vert_{L^1}\lesssim \sum_{B\in V_n}\Vert \psi_B\Vert_{L^1}
|Ph(B)| \lesssim \sum_{B\in V_n} |B| |Ph(B)|\lesssim \Vert h\Vert_{L^1}.
\end{align*}
If we apply this to $h=f-g$, then for some large $n\in \N$ we have
\begin{align*} \Vert \trace (Pf)-f\Vert_{L^1}\le
\Vert \trace (Pf)&-T_n(Pf)\Vert_{L^1}+\Vert T_n(P(f-g))\Vert_{L^1}+\\  &
\Vert T_n(Pg)-g\Vert_{L^1} +\Vert g-f\Vert_{L^1}
\lesssim \eps.
\end{align*}
Since $\eps>0$ was arbitrary, we conclude that $f=\trace (Pf)$ (in $L^1(Z)$, i.e., the functions agree almost everywhere on $Z$). 

For the second part we note that 
$$\frac{1}{|Z|}\int_Zf\, d\mu=T_0(Pf), $$
because $V_0=\{Z\}$. 
So if we use Lemma~\ref{le:trace}, then we  obtain 
\begin{align*}
 \biggl\Vert f-\frac{1}{|Z|}\int_Zf\,d\mu\biggr\Vert_{L^1}&=\Vert f-T_0(Pf)\Vert_{L^1}\\
 &\leq \sum_{n=0}^\infty
 \Vert T_{n+1}(Pf)-T_n(Pf)\Vert_{L^1}\\
 &\lesssim \|d(Pf) \|_{\ell^{p,\infty}}. \qedhere
\end{align*}
\end{proof}

\begin{proof}[Proof of Proposition~\ref{co:subspace}] Let  $f\in L^1(Z)$. Since $X$ is connected, we  have  $d(Pf)=0$ if
and only if $Pf$ is constant on $V$. This happens precisely  if $f$ is constant (almost everywhere on $Z$). In particular, $\| f\|_{A^p}=\|d(Pf)\|_{\ell^{p,\infty}}$ defines a norm on $A^p(Z)$, and the  linear map $f\mapsto d(Pf)$ gives an isometric embedding of  $A^p(Z)$ into $\ell^{p,\infty}(E).$
To finish the proof, it is enough to show that every Cauchy sequence $\{f_n\}$
  in 
$A^p(Z)$ converges to an element in $A^p(Z)$ (then the isometric  image of 
$A^p(Z)$ in $\ell^{p,\infty}(E)$ is closed).  Since $f_n$ is only well-defined up to a constant, we may assume that $\int_Zf_n\, d\mu=0$. Then 
Lemma~\ref{lem:frecov} shows that $$ \Vert f_n-f_k\Vert_{L^1} \lesssim  \Vert d(P(f_n-f_k))\Vert_{\ell^{p,\infty}}$$ for $k,n\in \N$. Hence 
$\{f_n\}$ is a Cauchy sequence in $L^1(Z)$, and so there exists $f\in L^1(Z)$ such that 
$ \Vert f_n-f\Vert_{L^1}\to 0$. Then  $Pf_n(B)\to Pf(B)$ for each $B\in V$, and so 
$d(Pf_n)(e)\to d(Pf)(e)$ for each $e\in E$.  Since $d(Pf_n)$ is uniformly bounded in 
$\ell^{p,\infty}(E)$,  this implies that $d(Pf)\in \ell^{p,\infty}(E)$ by Lemma~\ref{le:weakconv}, and so $f\in A^p(Z)$. 
Moreover, using the same lemma, we also see that  
$$\Vert d(Pf_n)-d(Pf)\Vert_{\ell^{p,\infty}}\to 0,$$ and so 
$f_n\to f$ in $A^p(Z)$  for $n\to \infty$ as desired. 
\end{proof}

 Given   $u\: V \ra \R$ we define  a type of maximal function $\mathcal{M}u\: V \ra [0,\infty]$
 by setting 
\beqla{defmax}
 (\mathcal{M}u)(B)=  
 \sum_{B'\in\bigcup_{k=n}^\infty V_k,\; 8B'\cap 8B\not=\emptyset}\frac{|B'|}{|B|}|u(B')|
\eeq
for $B\in V_n$, $n\in \N_0$. There is nothing special about the constant $8$ in the condition $8B'\cap 8B\not=\emptyset$. Its usefulness will become apparent later in  the proof of Corollary~\ref{cor:indfill}.

\begin{lemma}\label{lem:bddwk} Let $p>1$. Then there exists a constant $C\ge 0$ such that 
$$ \Vert \mathcal{M}u\Vert_{\ell^{p,\infty}}\le C \Vert u\Vert_{\ell^{p, \infty}}$$
for each $u \in \ell^{p,\infty}(V)$. 
\end{lemma} 
In other words, the sublinear operator $u\mapsto \mathcal{M}u$ is bounded in 
$\ell^{p,\infty}(V)$.

For $Z=\R^n$ this  lemma  is  essentially \cite[Lemma~4.6, p.~711]{ABH}; see also \cite[p.~269]{RS}.

\begin{proof}  First note that the operator $\mathcal{M}$ can be obtained 
from a kernel   $k\: V \times V\ra [0,\infty)$.  Indeed, we define  
$$\mathcal{K}(B,B')=|B'|/|B|$$ if $B\in V_n$, $B'\in V_k$, $k,n\in \N$, $k\ge n$, and $8B\cap 8B'\ne \emptyset$, and $\mathcal{K}(B,B')=0$ otherwise.  Then for $u\: V\ra \R$ and $B\in V$ we have 
$$\mathcal{M}u(B)=\sum_{B'\in V} \mathcal{K}(B,B') |u(B')|. $$
The main point is that this kernel $\mathcal{K}$ satisfies  a Schur condition (see \cite[pp.~86--87]{Wo}):  if $q>1$, $\alpha\in (0, 1/q)$ and $g(B)=|B|^\alpha$ for $B\in V$, then
\begin{align*} \sum_{B'\in V}\mathcal{K}(B,B')g(B')^{q/(q-1)}&=\frac{1}{|B|}\sum_{k=n}^\infty
\sum_{B'\in V_k, \, 8B\cap 8B'\ne \emptyset} |B'|^{1+\alpha q/(q-1)} \\
&\lesssim \sum_{k=n}^\infty 2^{-kQ\alpha q/(q-1)}\\ &\simeq 2^{-nQ\alpha q/(q-1)}\simeq  g(B)^{q/(q-1)}
\end{align*} 
for $B\in V_n$, $n\in \N_0$, 
and 
\begin{align*} \sum_{B\in V}g(B)^{q}\mathcal{K}(B,B')&=|B'|\sum_{n=0}^k \sum_{B\in V_n, \, 8B\cap 8B'\ne \emptyset} |B|^{-1+\alpha q}\\
&\simeq  |B'|\sum_{n=0}^k  2^{-nQ(-1+\alpha q)}\simeq  |B'| 2^{-kQ(-1+\alpha q)} \\
&\simeq 2^{-kQ\alpha q} \simeq |B'|^{\alpha q}\simeq g(B')^q
\end{align*} 
for $B'\in V_k$, $k\in \N_0$.  
This implies that the operator $\mathcal{M}$ is bounded in $\ell^q(V)$ for each $q>1$ and hence bounded in $\ell^{p,\infty}(V)$ by interpolation \cite[Chapter V.3]{SW}. 
\end{proof}

\begin{lemma}\label{pr:modout} Let $p>1$. 
\begin{itemize}
\item [\textnormal{(i)}] If  $u,u'\:V \ra \R$ with 
$du\in \ell^{p,\infty}(E)$ and $u-u'\in \ell^{p,\infty}(V)$,  then $du'\in \ell^{p,\infty}(E)$ 
and $\trace u =\trace u'.$

\smallskip 
\item [\textnormal{(ii)}] If $u\:V \ra \R$ and 
 $du\in \ell^{p,\infty}(E)$,  then
\begin{align*}
\|u-P(\trace u)\|_{\ell^{p,\infty}}&\lesssim \|du \|_{\ell^{p,\infty}}\quad \text{ and }\\
\|d(P(\trace u)\|_{\ell^{p,\infty}}&\lesssim \|du \|_{\ell^{p,\infty}}. 
\end{align*}
In particular,  $u-P(\trace u)\in \ell^{p,\infty}(V)$.

\end{itemize}
\end{lemma}
\begin{proof} (i) We have 
\begin{align*}
 \Vert du'\Vert_{\ell^{p,\infty}}
&\le  \Vert d(u'-u)\Vert_{\ell^{p,\infty}}+ \Vert du\Vert_{\ell^{p,\infty}} \\ &\lesssim 
 \Vert u'-u\Vert_{\ell^{p,\infty}}+ \Vert du\Vert_{\ell^{p,\infty}}<\infty,
\end{align*}
and so both $\trace u$ and  $\trace u'$ are defined by Lemma~\ref{le:trace}.

If   $w\: V \ra \R$ is a function with $w\in \ell^{p,\infty}(V)$, then
 $|w(B)|$ is  small if the level of $B\in V$ is large enough. 
Hence  $\sup_{x\in Z}|(T_nw)(x)|\to 0$ and so
$\| T_nw\|_{L^1}\to 0$ as $n\to\infty$. If we apply this to $w=u-u'$, the  claim follows.

\smallskip
(ii) If $u$ is as in the statement, then an estimate as in \eqref{diffest} implies that 
for $B\in V_n$ with $n\in \N_0$ we have,  
\begin{align*}
|u(B)&-P(\trace u) (B)| \leq  \\
&\biggl|u(B)-\frac {1}{|B|}\int_BT_nu\, d\mu\biggr|+\sum_{k=n}^\infty \frac{1}{|B|}
\int_B|T_ku-T_{k+1}u|\, d\mu\\
&\leq \sum_{B'\in\bigcup_{k=n}^\infty V_k,\; B'\cap B\not=\emptyset}\frac{|B'|}{|B|}\tilde du(B')
\le\mathcal{M}(\tilde du )(B).\nonumber
\end{align*}
Here $\mathcal{M}(\tilde du)$ is as defined in \eqref{defmax}. The first inequality  now follows from 
Lemma~\ref{lem:bddwk} and \eqref{dduin}.  

For the second inequality note that 
\begin{align*} \Vert d(P(\trace u))\Vert_{\ell^{\infty, p}}
&\lesssim  \Vert d (u-P(\trace u))\Vert_{\ell^{\infty, p}}+ \Vert du\Vert_{\ell^{\infty, p}}\\
&\lesssim  \Vert u-P(\trace u)\Vert_{\ell^{\infty, p}} + \Vert du\Vert_{\ell^{\infty, p}}\\
&\lesssim  \Vert du\Vert_{\ell^{\infty, p}}. \qedhere
\end{align*}
\end{proof}

\begin{proof}[Proof of Theorem~\ref{main2}] If $[u]\in \Wk^{1,p}(X)$, then 
$du\in \ell^{p,\infty}(E)$, and so $\trace u\in L^1(Z)$ is defined by Lem\-ma~\ref{le:trace}. Moreover, by Lem\-ma~\ref{pr:modout}~(ii) we have 
\beqla{eq:nobdd} \Vert d(P(\trace u))\Vert_{\ell^{\infty, p}}
\lesssim  \Vert du\Vert_{\ell^{\infty, p}}<\infty, 
\eeq
and so $\trace u\in A^p(Z)$. Note that $\trace u$   is well-defined in $A^p(Z)$ and  depends only  on 
$[u]$; indeed, if $u'\in [u]$, then there exists a constant $c\in \R$ such that $u-u'-c\in 
\ell^{p,\infty}(V)$. Then by Lemma~\ref{pr:modout}~(i) we have  $\trace u=\trace (u'+c)=c+\trace u'$.  So we see  that 
$[u]\in \Wk^{1,Q}(X)\mapsto \trace u\in A^p(Z)$ defines a linear map $\trace \:  
\Wk^{1,p}(X)\ra A^p(Z)$. 

It is clear that $[u]\in  \Wk^{1,p}(X)\mapsto \Vert u\Vert_{\Wk^{1,p}} $ defines a semi-norm  on 
 $\Wk^{1,p}(X)$. Suppose  that $ \Vert u\Vert_{\Wk^{1,p}}  =0$ for $[u]\in  \Wk^{1,p}(X)$. By taking the infimum in \eqref{eq:nobdd} over all representatives in $[u]$ we see that 
 $$\Vert \trace u\Vert_{A^p(Z)}=\Vert d(P(\trace u))\Vert_{\ell^{\infty, p}}=0, $$ and so $\trace u$ is constant; but then $u\in \ell^{p,\infty}(V)+\R$ by  
  Lemma~\ref{pr:modout}~(ii), showing that $[u]=0$. Hence  $[u]\mapsto \Vert u\Vert_{\Wk^{1,p}}  $ is a norm on $\Wk^{1,p}(X)$. 
  Inequality~\eqref{eq:nobdd} also shows that 
  $\trace \: \!\!
\Wk^{1,p}(X)\ra A^p(Z)$ is a bounded linear operator from the normed space 
 $\Wk^{1,p}(X)$ into the Banach space $A^p(Z)$.

It is immediate from Lemma~\ref{lem:frecov} that the map $f\in A^p(Z)\mapsto [Pf]\in \Wk^{1,p}(X)$ is a bounded linear operator $\mathcal{P}\:A^p(Z) \ra \Wk^{1,p}(X)$ such that the composition 
$\trace \circ \mathcal{P} $ is the identity on $A^p(Z)$. If $[u]\in 
 \Wk^{1,p}(X)$, then, by what we have seen,  $d(P(\trace u)) \in \ell^{p,\infty}(E)$. Moreover, 
 $P(\trace u)-u\in \ell^{p,\infty}(V)$ by Lemma~\ref{pr:modout}~(ii). Hence $[u]=[P(\trace u)]$ which shows that $\mathcal{P}\circ \trace$ is the identity on $\Wk^{1,p}(X)$. It follows that $\trace$ is an isomorphism between the normed spaces $\Wk^{1,p}(X)$ and $A^p(Z)$ with the inverse given by the operator $\mathcal{P}$. Since $A^p(Z)$ is a  Banach space,  $\Wk^{1,p}(X)$ is a Banach space as well. The proof is complete. \end {proof}
 
 In order to show that $A^p(Z)$ does not depend on the filling $X$ in any essential
 way, we need a lemma. 
 
\begin{lemma}\label{lem:strweaest} Let $p>1$, and   $f\in L^1(Z)$ with 
$d(Pf)\in \ell^{p,\infty}(E)$. Define $Df\: V\ra \R$ as 
$$ (Df)(B)=\frac{1}{|8B|}\int_{8B}|f-f_B|\,d\mu$$
for $B\in V$. Then 
$$ \Vert Df\Vert_{\ell^{p,\infty}}\le C\Vert d(Pf) \Vert_{\ell^{p,\infty}}, $$
where $C\ge 0$ is constant independent of $f$.
\end{lemma}

\begin{proof} Let $u=Pf$. Then by Lemma~\ref{le:trace} and Lemma~\ref{lem:frecov} we have  $f=\trace (Pf)=\lim_{n\to \infty}T_nu$ with convergence in $L^1(Z)$ and pointwise almost everywhere on $Z$. It follows that 
for $B\in V_n$, $n\in \N_0$, we have 
\begin{align*}
 (Df)(B)&=\frac{1}{|8B|}\int_{8B}|f-f_B|\,d\mu\\
 &\le \frac{1}{|8B|} \int_{8B}|u(B)-T_nu|\, d\mu+ \frac{1}{|8B|}\sum_{k=n}^\infty
 \int_{8B}
 |T_{k+1}u-T_ku|\, d\mu \\
& \lesssim \frac{1}{|B|} \sum_{B'\in\bigcup_{k=n}^\infty V_k,\; B'\cap 8B\not=\emptyset}|B'| \tilde d u(B') \le \mathcal{M}(\tilde du)(B),
\end{align*}
where $\mathcal M$ is defined as in \eqref{defmax}. By Lemma~\ref{lem:bddwk} and \eqref{dduin} this implies 
$$ 
 \Vert Df\Vert_{\ell^{p,\infty}}\lesssim \Vert  \mathcal{M}(\tilde du) \Vert_{\ell^{p,\infty}}  
 \lesssim \Vert \tilde du \Vert_{\ell^{p,\infty}}\lesssim \Vert du \Vert_{\ell^{p,\infty}}, 
$$
and the claim follows. 
\end{proof}

\begin{corollary} \label{cor:indfill} For $p>1$ is space is $\dot A^p(Z)$ is independent of the hyperbolic filling
with comparability of semi-norms for different fillings. 
\end{corollary}

\begin{proof} Suppose  $ X'=(V', E')$ is a hyperbolic filling of $Z$ constructed  as in the beginning of this section using  possibly different  maximal $2^{-n}$-separated sets $Z_n'$.  Then again $V'$ is  a collection of balls, and we denote by $P'f\: V'\ra \R$ the Poisson extension of a function $f\in L^1(Z)$, and by $du'\: E'\ra \R$ the gradient of a function $u'\: V'\ra \R$. 
 It  suffices to show that 
 $\Vert  d(P'f) \Vert_{\ell^{p,\infty}} \simeq \Vert d(Pf) \Vert_{\ell^{p,\infty}} $
 for $f\in L^1(Z)$. For this it is enough to show that 
 $\Vert  d(P'f) \Vert_{\ell^{p,\infty}} \lesssim  \Vert d(Pf) \Vert_{\ell^{p,\infty}} $
 by symmetry. 
 
 So suppose $f\in L^1(Z)$, and define $u'=P'f$. Let $e'\in E'$ be arbitrary,
 and $n\in \N_0$ be the level of $Q:=e'_+$. Then there exists $B\in V_n$ such that 
 $B\cap Q\ne \emptyset. $ Since the radii of $B$ and $Q$ agree, and the radius of 
 $Q'=e'_-$ is at most twice as large as the radius of $Q$ we have 
 $ Q,Q'\sub 8B$ (this inclusion is the reason why we choose the constant $8$ in
 the definition  \eqref{defmax} of $\mathcal{M}$). 
 It follows that 
 \begin{align} 
  \label{xy} |d(P'f)(e')|  &  =|f_Q-f_{Q'}| \le |f_Q-f_B|+ |f_{Q'}-f_B|   \\  
 &\le \frac{1}{|Q|} \int_Q |f-f_B|\,d\mu+ \frac{1}{|Q'|} \int_{Q'}|f-f_B|\,d\mu \nn  \\ 
  &\lesssim \frac{1}{|8B|}\int_{8B}|f-f_B|\,d\mu =(Df)(B). \nn
  \end{align}
   Here we used that $|Q|\simeq |Q'|\simeq |8B|$. 
Note that for a given ball $B\in V$ there  is only a uniformly bounded number of edges $e'\in E'$ so that $B$ can appear in \eqref{xy}. Hence Lemma~\ref{le:weakbdd} and 
Lemma~\ref{lem:strweaest} imply that
$$ \Vert  d(P'f) \Vert_{\ell^{p,\infty}} \lesssim     \Vert  Df \Vert_{\ell^{p,\infty}}\lesssim\Vert d(Pf) \Vert_{\ell^{p,\infty}}  $$ as desired.
\end{proof}

\begin{remark}\label{rem:fill}
\textnormal{In the previous proof  we insisted on all hyperbolic fillings to be constructed 
from the scale range $2^{-n}$.  There is no difficulty in extending the scope of 
Corollary~\ref{cor:indfill} to hyperbolic fillings obtained from  an arbitrary scale 
range $\Lambda^{-n}$ with $\Lambda>1$. This only requires   some adjustment 
of the constant $8$ in the definition  \eqref{defmax} of the maximal operator $\mathcal{M}$.  This leads  to more general versions of  Lemma~\ref{lem:strweaest}  that  
can be used to extend Corollary~\ref{cor:indfill}. We omit the details, because this only leads to technicalities adding nothing of substance.
}\end{remark}

 \section{Sobolev spaces}\label{s:PI}
 In this section  $Z=(Z,d,\mu)$ is again   an  Ahlfors $Q$-regular compact metric measure space with $Q>0$, and $X=(V,E)$  the hyperbolic filling of $Z$.
 
If   $\alpha>0$ and $f\in L^1(Z)$, then  we denote by $ \mathcal{D}_\alpha(f)$ the  set of all measurable functions $g\ge 0$ on $Z$ for which  there exists a set $N\sub Z$ with $|N|=0$ such that 
\begin{equation}\label{def:Haj} |f(x)-f(y)|\le d(x,y)^\alpha(g(x)+g(y))
\end{equation} 
for  all $x,y\in Z\setminus N$. 

For $\alpha>0$ and $p\ge 1$ we denote by $\dot M^{\alpha,p}(Z)$ the 
   (homogeneous) {\em fractional Haj\l asz-Sobolev space}  consisting  of all integrable functions
$f\in L^1(Z)$ such that there exists a function $g\in  \mathcal{D}_\alpha(f)$ with  $g\in L^p(Z)$; in other words, $f\in \dot M^{\alpha,p}(Z)$ if there exists a non-negative function $g\in L^p(Z)$ such that  \eqref{def:Haj} is true for almost every $x$ and $y$ in $Z$.  
The spaces $\dot M^{\alpha,p}(Z)$ were introduced in \cite{Ha} for $ \alpha=1$ and in \cite{Ya} for arbitrary $\alpha>0$.
We define a semi-norm on $\dot M^{\alpha,p}(Z)$ by setting
$$\Vert f\Vert_{\dot M^{\alpha,p}}=\inf \{\Vert g\Vert_{L^p}: g\in \mathcal{D}_\alpha(f)\}.$$

\begin{lemma}\label{lem:convres} Let $p>1$, $\alpha>0$,  and  $f\in L^1(Z)$. Suppose 
 $\{f_n\}$ is a sequence in  $L^{1}(Z)$ 
such that   $f_n(x)\to f(x)$ as $n\to \infty$ for almost every $x\in Z$. 
Then $$\Vert f\Vert_{\dot M^{\alpha,p}}\le \liminf_{n\to \infty} 
\Vert f_n\Vert_{\dot M^{\alpha,p}}.$$ 
In particular, if the sequence $\{f_n\}$ lies in $\dot M^{\alpha,p}(Z)$ and 
$\Vert f_n\Vert_{\dot M^{\alpha,p}}$ is uniformly bounded for $n\in \N$, then 
$f\in \dot M^{\alpha,p}(Z)$.
\end{lemma}

\begin{proof} If  $C_0:=\liminf_{n\to \infty} 
\Vert f_n\Vert_{\dot M^{\alpha,p}}=\infty$ there is nothing to prove.
So let us assume that $C_0<\infty $ and  pick $C>C_0$. Then by passing to a subsequence if necessary, we may assume that $\Vert f_n\Vert_{\dot M^{\alpha,p}}< C$ for all
$n\in \N$. Hence for each $n\in \N$ there exists $g_n\in \mathcal{D}_\alpha(f_n)$ such that $\Vert g_n\Vert_{L^p}<C$.  So the sequence $\{g_n\}$ is uniformly bounded in $L^p(Z)$. Again by passing to a subsequence, we may assume that $g_n\to g\in L^p(Z)$ with respect to the weak-$\ast$ topology in $L^p(Z)$. By Mazur's lemma 
\cite[Section V.1]{Yo} we can find 
 a sequence  $\{\tilde g_n\}$ in $L^p(Z)$, where each function $\tilde g_n$ is a finite convex combination of functions $g_k$ with $k\ge n$,  such that we have norm-convergence 
$\tilde g_n \to g$ in $L^p(Z)$. Then $\Vert \tilde g_n\Vert_{L^p}\le C$ for all $n\in \N$ and hence also $\Vert g\Vert_{L^p}\le C$. 

Note  that condition \eqref{def:Haj} passes to convex combinations of  pairs $(f,g)$. So by using   the same coefficients  as in the definition of $\tilde g_n$, we can find a convex combination $\tilde f_n$ involving the functions 
$f_k$ with $k\ge n$ such that $\tilde g_n\in \mathcal{D}_\alpha(\tilde f_n)$. Then  we still have $\tilde f_n(x)\to f(x)$ as $n\to \infty$ for almost every $x\in Z$.  By passing to yet another subsequence, we may assume that $\tilde g_n(x)\to g(x)$ for almost every $x\in Z$. 
Then by using \eqref{def:Haj} for the pair $(\tilde f_n, \tilde g_n)$ and passing to the pointwise limit in $x$ and $y$ outside a suitable set $N\sub Z$ with $|N|=0$, we see that 
$g\in \mathcal{D}_\alpha(f)$. So $\Vert f\Vert_{\dot M^{\alpha,p}}\le C$.  The claim follows by letting $C\to C_0$.  
\end{proof}

\begin{lemma}\label{conslem} Let $p\ge 1$, $\alpha>0$,  and  $f\in L^1(Z)$.
If $\Vert f\Vert_{\dot M^{\alpha,p}}=0$, then $f$ is constant almost everywhere on $Z$.
\end{lemma}

\begin{proof} If $\Vert f\Vert_{\dot M^{\alpha,p}}=0$, then there exists a sequence $\{g_n\}$ of functions in $\mathcal{D}_\alpha(f)$ with $\Vert g_n\Vert_{L^p}\to 0$ as $n\to \infty$. By passing to a suitable subsequence, we may assume that $g_n(x)\to 0$ as $n\to \infty$ for almost every $x\in Z$. A limiting argument then implies that the constant function $g_\infty\equiv 0$ lies in $\mathcal{D}_\alpha(f)$. Hence  $f$ is constant almost everywhere on $Z$. 
\end{proof}

If $g\in L^1(Z)$, then we  we define the  {\em Hardy-Littlewood maximal function} by 
\begin{equation}\label{HKMAX} (Mg)(x)= \sup_{ x\in B} \frac1{|B|}\int_{B}|g| \,d\mu,
\end{equation}
for $x\in Z$, where the supremum is taken over all balls $B\subset Z$ with $x\in B$. 
It is a standard fact that for $p>1$ the operator $f\mapsto M(f)$ is bounded in $L^p(Z)$
\cite[pp.~10--12]{He}. 

\begin{proposition}\label{Hajincl} Let $p>1$, and $\alpha=Q/p$. Then $\dot M^{\alpha,p}(Z)\subset \dot A^p(Z)$. Moreover, there exists a constant $C\ge 0$ such that 
$$ \Vert f \Vert_{\dot A^p(Z)}\le  C\Vert  f \Vert_{\dot M^{\alpha,p}}$$ for all $f\in \dot M^{\alpha,p}(Z)$.
\end{proposition}

\begin{proof} Let $p>1$, $\alpha=Q/p$,  and $f\in \dot M^{\alpha,p}(Z)$. Then there exists 
 $g\in \mathcal{D}_\alpha(f)$ with $g\in L^p(Z)$.
 If  $u=Pf$, then we have  to show that $du\in \ell^{p,\infty}(E)$ with a  suitable norm bound.  
  
For   $\lambda >0$ let   us denote by $E(\lambda)$ the collection of all edges  $e\in E$ such that $|du(e)|>\lambda$ and by $V(\lambda)$ the collection of all corresponding balls  $B=e_+$, where $e\in E(\lambda)$. Then we have 
\beqla{eq:counting}
\# E(\lambda)\lesssim \#B(\lambda) \le \int_Z\Big( \sum_{B\in V(\lambda)} \frac1{|B|} \chi_{B}\Big)\,d\mu,
\eeq
where $\chi_B$ denotes the characteristic function of $B$. 
Moreover, by the  dyadic structure of the balls  $B$  we have  
\begin{equation} \label{ballest} \sum_{B\in V(\lambda)} \frac1{|B|  }  \chi_{B}(z) \lesssim \frac {1}{|B_z|}, 
\end{equation}
where $B_z$ is a  ball  in $V_\lambda$ of smallest radius that contains $z\in Z$ if such a ball exists. 

We consider $e\in E$, and let $B=e_+$ and $B'=e_-$. Note that $\diam(B)\simeq |B|^{1/Q}$ and $\diam(B)^\alpha\simeq |B|^{1/p}$. Then for a sufficiently large constant  $\Lambda\ge 1$ independent of $e$ we have 
\begin{align*}|du(e)| &=|u(B)-u(B')| =| f_B-f_{B'}| \\
&\le \frac{1}{|B|\cdot |B'|} \int_B\int_{B'} |f(x)-f(y)|\,d\mu(x)d\mu(y)\\
&\lesssim  \frac{\diam(B)^\alpha}{|B|\cdot |B'|} \int_B\int_{B'}(g(x)+g(y))\,d\mu(x)d\mu(y)\\
&\lesssim |B|^{1/p}  \frac1{|\Lambda B|}\int_{\Lambda B}g\,d\mu\lesssim
|B|^{1/p} (Mg)(z),
\end{align*}
whenever $z\in B$. Here $Mg$ is the Hardy-Littlewood maximal function of $g$ as in \eqref{HKMAX}. 
The previous inequality shows  that  
\begin{equation}\label{eq:maximal1}
|du(e)|>\lambda \Rightarrow \frac1{|B|}\lesssim \frac{1}{\lambda^p} (M(g)(z))^{p} \text{ for all  $z\in B=e_+$.} 
\end{equation}
Note that this implies that if $z\in Z$ and there exist balls $B\in V(\lambda)$ of arbitrarily small radius with $z\in B$, then $M(g)(z)=+\infty$. Otherwise, we can apply \refeq{ballest}. This leads to the inequality 
$$ \sum_{B\in V(\lambda)}\frac1{|B|  }  \chi_{B}\lesssim \frac{1}{\lambda^p}  M(g)^{p}. $$ 
Substituting this into \refeq{eq:counting} and integrating over $Z$,  we  finally obtain
\beqla{eq:normbdd1}
\# E(\lambda)\lesssim  \frac{1}{\lambda^p}   \int_X M(g)^{p}\, d\mu \lesssim
 \frac{1}{\lambda^p}   \int_X  g^{p}\,d\mu. 
\eeq
Here we used $p>1$ and that  the   maximal function is bounded in $L^p(Z)$. Inequality \eqref{eq:normbdd1}
 implies 
$$\| f\|_{\dot A^p}=\Vert d(Pf)\Vert_{\ell^{p,\infty}} \lesssim \Vert g\Vert_{L^p}<\infty,$$
and so $f\in A^p(Z)$.  Taking the infimum over all $g\in \mathcal{D}_{\alpha}(f)$ in the last inequality, we conclude that  
$\| f\|_{A^p}\lesssim\| f\|_{\dot M^{\alpha,p}}$ as desired. 
\end{proof}

 For a Lipschitz function 
$f\: Z\ra \R$ and $x\in Z$ we define 
$$ \lip_xf=\limsup_{r\to 0^+} \sup_{y\in  B(x,r)} |f(y)-f(x)|/r, 
$$
and a measurable function $ \lip\,f\: Z\ra \R$ by setting 
$ (\lip\,f)(x)=\lip_xf$ for $x\in Z$. 
For $p\ge 1$ we say that  $Z$ supports a {\em $p$-Poincar\'e inequality} if there exist constants $C\ge 0$ and $\Lambda\ge 1$ such that 
$$ \frac 1{|B|}\int |f-f_B|\, d\mu\le C R \biggl( \frac 1{|\Lambda B|} \int_{\Lambda B} 
 (\lip\, f)^p\, d\mu\biggr)^{1/p} $$
for every Lipschitz function $f\: Z \ra \R$ and 
every  ball $B\subset Z$, where $R$ is the radius of $B$. Note that our definition is equivalent to the more standard one given in \cite{He}; see  
  \cite[Theorem~2] {Ke}.
  
  The following lemma is essentially well-known.

\begin{lemma} \label{lem:LIP}
Let $(Z,d,\mu)$ be an Ahlfors $Q$-regular compact metric measure space that supports a $Q$-Poincar\'e inequality, where $Q> 1$. Then there exists a constant $C>0$ 
such that 
$$  \Vert f\Vert_{\dot M^{1,Q}} \le  C  \Vert \! \lip f\Vert_{L^Q}$$
for every Lipschitz function  $f\: Z \ra \R$.  \end{lemma}

\begin{proof} By \cite{KZ} we know that there exists $r\in [1,Q)$ such that 
$(Z,d,\mu)$ supports an $r$-Poincar\'e inequality. 

Let $f\: Z \ra \R$ be a Lipschitz function and $x,y\in Z$ with $x\ne y$ be arbitrary. Define $\rho=\lip f$ and 
 $B_n=B(x, 2^{-n} d(x,y))$ for $n\in \N_0$. Then the $r$-Poincar\'e inequality 
 implies that for some constant $\Lambda\ge 1$ we have
 \begin{align*} |f(x)-f_{B_0}|  & \le \sum_{n=0}^\infty |f_{B_n}-f_{B_{n+1}}| \\
& \lesssim d(x,y) \sum_{n=0}^\infty 2^{-n}  \biggl( \frac {1} {|\Lambda B_n|} \int_{\Lambda B_n} 
(\lip f)^r\, d\mu  \biggr)^{1/r}\\
&\lesssim d(x,y) [M(\rho^r)(x)]^{1/r}. 
\end{align*}
Here $M$ denotes the Hardy-Littlewood maximal function defined as in \eqref{HKMAX}. 

If we define $B_0'=B(y, d(x,y))$, then very similar estimates lead to 
$$  |f(y)-f_{B'_0}| \lesssim  d(x,y)  [M(\rho^r)(y)]^{1/r}$$ 
and 
$$ |f_{B_0}-f_{B'_0}| \lesssim  d(x,y)  [M(\rho^r)(x)]^{1/r}. $$
Hence 
$$ |f(x)-f(y)|\lesssim d(x,y) \big(  [M(\rho^r)(x)]^{1/r}+ [M(\rho^r)(y)]^{1/r}\big),  $$
where the implicit multiplicative constant is independent of $f$, $x$, $y$.
It follows that for some constant $k_0$ independent of $f$ we have 
 $k_0[M(\rho^r)]^{1/r} \in \mathcal{D}_1(f)$. 
 Hence 
 \begin{align*} \Vert f\Vert_{\dot M^{1,Q}}^Q &\lesssim \Vert [M(\rho^r)]^{1/r}\Vert_{L^Q}^Q
 =\int_Z [M(\rho^r)]^{Q/r}\, d\mu\\
 &\lesssim \int_Z[\rho^r]^{Q/r}\, d\mu= \int_Z(\lip f)^{Q}\, d\mu=
 \Vert \! \lip f\Vert^Q_{L^Q}.
  \end{align*}
Here we used that $Q/r>1$ which implies that $g\mapsto M(g)$ is bounded in $L^{Q/r}(Z)$. The claim follows. 
\end{proof} 

The previous argument actually shows that if $p\le Q$ is sufficiently close to $Q$ (namely if $p\in (r,Q]$),  then 
$\Vert f\Vert_{\dot M^{1,p}} \lesssim \Vert \! \lip f\Vert_{L^p}$.

\begin{proof}[Proof of Theorem~\ref{main}] Let  $Z$ be as in the statement, and  $X=(V,E)$ be the hyperbolic filling of $Z$ as discussed in Section~\ref{s:Ahl}.
 By Proposition~\ref{Hajincl} we know that $\dot M^{1, Q}(Z)\sub \dot A^{Q}(Z)$ with a suitable semi-norm bound.

For  the other direction, let us assume that $f\in L^1(Z)$ and  $u=Pf$ satisfies $du\in \ell^{Q,\infty}(E)$. Fix $N\in \N$ and consider the  gradient functions  $\tilde d_nu$ on levels  
$n\in [N,2N]$.  The total number of   vertices  corresponding to these levels is $\simeq 2^{2NQ}$. 
Hence 
$$ \sum_{n=N}^{2N} \Vert \tilde d_nu\Vert_{\ell^Q}^Q\lesssim N \Vert  du\Vert^Q_{\ell^{Q,\infty}} $$
as follows from  \eqref{dduin} and the first inequality in  Lemma~\ref{le:weaktostrong}~(ii). In particular, there exists $n\in [N,2N]$ such that 
$$\Vert \tilde d_{n}u\Vert_{\ell^Q}\lesssim  \Vert du\Vert_{\ell^{Q,\infty}}.$$ By using this fact on a non-overlapping sequence of such intervals  $[N,2N]$ one can find a sequence $\{n_k\}$ in $\N$ with
$n_k\to \infty$ and 
$$\Vert \tilde d_{n_k}u\Vert_{\ell^Q}\lesssim  \Vert du\Vert_{\ell^{Q,\infty}}$$
for all $k\in \N$. Then  $\Vert  \tilde d_{n_k}u\Vert_{\ell^Q}$ stays uniformly bounded as $k\to \infty$.

By Lemma~\ref{le:trace} and Lemma~\ref{lem:frecov} we have  
$$f=\trace (Pf) =\trace u=
\lim_{k\to\infty}T_{n_k}u, $$ where the limit is in $L^1(Z)$ and pointwise almost everywhere in $Z$.  So  $f\in \dot M^{1,Q}(Z)$  will follow from Lemma~\ref{lem:convres}, if we can show   that the norms
$\|T_{n_k}u\|_{\dot M^{1,Q}}$ are uniformly bounded  for $k\in \N$. 

To see this,  let $n\in\N$ and $x\in Z$. We pick $B\in V_n$ with $x\in B$. Then for 
$y\in Z$ we have,   
$$ |T_nu(y)-T_nu(x)| \le \sum_{B'\in V_n}  |u(B')-u(B)|\cdot  |\psi_{B'}(y)-\psi_{B'}(x)|, $$
and so 
\begin{align}\label{lipbdd} 
\lip_x(T_{n}u)&= \limsup_{r\to 0^+}\, \sup\{ |T_nu(y)-T_nu(x)|/r: y\in Z, \, d(y,x)< r\}
\nonumber \\
&\le \sum_{B'\in V_n,\, x\in B'}  |u(B')-u(B)| \lip_x(\psi_{B'})  \\
&\lesssim 2^{n} \sum_{B'\in V_n,\, x\in B'}  |u(B')-u(B)|\lesssim  
2^{n}   \tilde du(B).\nonumber 
\end{align}
Here we used \eqref{suppBincl} and \eqref{Lipbdd}. 
Hence by Lemma~\ref{lem:LIP} we have 
$$\Vert T_nu\Vert_{\dot M^{1,Q}}^Q\lesssim  \Vert \! \lip(T_{n}u)\Vert_{L^Q}^Q \lesssim  \sum_{B\in V_n} \tilde du(B)^Q=
\Vert \tilde d_nu\Vert^Q_{\ell^{Q}}. 
$$
Applying this inequality on the subsequence $\{n_k\}$, we see that 
\begin{align*}
 \Vert f\Vert_{\dot M^{1,Q}} &\le \liminf_{k\to \infty} \Vert T_{n_k}u\Vert_{\dot M^{1,Q}}
 \lesssim \liminf_{k\to \infty} \Vert\! \lip(T_{n_k}u)\Vert_{L^Q}
 \\
& \lesssim  \liminf_{k\to \infty} \Vert \tilde d_{n_k}u\Vert_{\ell^{Q}}\lesssim \Vert du\Vert_{\ell^{Q,\infty}}< \infty.
\end{align*}
So $f\in \dot M^{1,Q}(Z)$, and we also  conclude that  $\| f\|_{\dot M^{1,Q}} \lesssim \| f\|_{\dot A^Q}$. The  proof is complete.
\end{proof}

\begin{proposition}\label{prop:consfunc}
Under the assumptions of Theorem~\ref{main}, the space $\dot A^p(Z)$  
 consists  only  of constant functions for $1< p<Q$. \end{proposition}
 
\begin{proof} Let $1< p<Q$,  $f\in A^p(Z)$, and $u=Pf$.  We pick $q\in (p,Q)$ close to $Q$. As in the proof of Theorem~\ref{main},  one shows that
there exists a sequence $\{n_k\}$ in $\N$ with $n_k\to \infty$ as $k\to \infty$ such that 
$\Vert \tilde d_{n_k}u\Vert_{\ell^p}\lesssim 1.$ Since $q>p$, we then  also have 
$\Vert \tilde d_{n_k}u\Vert_{\ell^q}\le \Vert \tilde d_{n_k}u\Vert_{\ell^p}\lesssim 1.$
 Inequality  \eqref{lipbdd}  implies  that 
 $$ \Vert\! \lip(T_nu)\Vert_{L^q}^q \lesssim 2^{-(Q-q)n}\Vert \tilde d_n u\Vert_{\ell^q}^q,
 $$
and so $\Vert  \!\lip(T_{n_k}u)\Vert_{L^q} \to 0$ as $k\to \infty$. If 
$q$ is sufficiently close to $Q$, then  we have (see the remark after the proof of 
Lemma~\ref{lem:LIP}),
\begin{align*} 
\Vert f\Vert_{\dot M^{1,q}}&\le \liminf_{n\to \infty} \Vert T_nu\Vert_{\dot M^{1,q}}
\lesssim \liminf_{n\to \infty} \Vert \!\lip(T_nu)\Vert_{L^q} \\
&\le \lim_{k\to \infty} \Vert \!\lip(T_{n_k}u)\Vert_{L^q}=0.
\end{align*} 
Hence $\Vert f\Vert_{\dot M^{1,q}}=0$, and so $f$ is a constant function by Lemma~\ref{conslem}.
\end{proof}

\section{The Euclidean setting}\label{s:rn}
We will now discuss results in the spirit of Theorem~\ref{main} in the standard Euclidean  setting. We choose the simplest framework, where the underlying space is $Z=\R^n$, $n\in \N$,    equipped  with Lebesgue measure. 

We denote Lebesgue measure of a measurable set $M\sub \R^n$ by $|M|$, and use $dx$ to indicate Lebesgue measure in integrals. 
We consider the  hyperbolic filling $X$ of $Z=\R^n$ as given by upper-halfspace $\R_+^{n+1}=\{(x,t): x\in \R^n,\, t>0\}$. In contrast to our earlier discussion, in this section  we consider the {\em (classical) Poisson extension } $u=Pf$ of a function $f\in L^p(\R^n)$, where $p\in [1,\infty].$ If 
 $$P_t(x)=c_n \frac {t}{(t^2+|x|^2)^{(n+1)/2}}, \quad x\in \R^n,\, t>0,$$ is the Poisson kernel, where $c_n=\Gamma((n+1)/2)/\pi^{(n+1)/2}$, then $u$ is given by 
the convolution 
 $u(x,t)=(P_t\ast f)(x)$ for $(x,t)\in \R_+^{n+1}$ (see \cite[Chapter~3.2]{St} for general background).  
 
If $\mu$ is a measure on  set $A\subset \R^n$ and $p\ge 1$, then we define the {\em weak-type space} $L^{p, \infty}(A, \mu)$ as the space of all measurable  
functions $f\: A \ra \R$ for which  there exists a constant $C\ge 0$ such that 
$$ \mu(\{ x\in A: |f(x)|>\lambda\})\le C^p/{\lambda^p} $$
for all $ \lambda>0$. We denote by $\Vert f\Vert_{L^{p,\infty}(A,\mu)}$ the infimum of all  constants $C\ge 0$ for which this inequality is valid.

We can now formulate a result, whose basic idea  is found \cite[Appendix]{CST}. It provides  a somewhat surprising concrete embedding of $L^p$-spaces into subspaces of  the  
weak-type spaces $L^{p,\infty}.$ 

\begin{proposition}\label{prop:fuweak}
Let  $p>1$ and $s>0$. If $f\in L^p (\R^n )$, then the Poisson extension $u=Pf$ satisfies
\beqla{eq:poissonwe}
t^{s/p}u(x,t)\in L^{p,\infty}\big (\R_+^{n+1},t^{-(s+1)}dx\,dt).
\eeq
Conversely, if a harmonic function $u$ on $\R_+^{n+1}$ satisfies $\refeq{eq:poissonwe}$, then $u=Pf$, where $f\in L^p (\R^n )$.
\end{proposition}

Of course, here $t^{s/p}u(x,t)$   stands for the function on  $\R_+^{n+1}$ given by
$$(x,t)\mapsto t^{s/p}u(x,t). $$ We will also  make a similar abuse of notation below. 

\begin{proof} We denote by $\mu_s$ the measure on $ \R_+^{n+1}$ given by 
$$d \mu_s(x,t)=t^{-(s+1)}\,dx\, dt.$$

Let $f\in L^p (\R^n)$. Then $Mf\in L^p(\R^n)$ for the Hardy-Littlewood maximal function of $f$, and a standard estimate for the Poisson extension \cite[p.~62]{St} gives     the bound 
\beqla{eq:maximal}
|u(x,t)|\le c  (Mf)(x)\eeq
for each $(x,t)\in \R_+^{n+1}$, where $c>0$ is a constant independent of $f$ and $x$.

Fix   $x\in \R^n$ and $\lambda>0$. If for $t>0$ we have 
$ |t^{s/p}u(x,t)|> \lambda $, then it follows that 
$$ t \ge   \biggl(\frac{\lambda}{c(Mf)(x)}\biggr)^{p/s}.
$$
We obtain
\begin{align*}
\mu_s \{(x,t)\in \R_+^{n+1}: |t^{s/p}&u(x,t)|> \lambda \}\\
&\leq \int_{\R^n}\left(\int^\infty_{(\lambda/c(Mf)(x))^{p/s}}\frac{dt}{t^{s+1}}\right)dx\\
& \leq \frac{c^p}{s\lambda^{p}}\int_{\R^n}[(Mf)(x)]^p\, dx
\lesssim  \frac1{\lambda^{p}}\|f \|_{L^p}^p,
\end{align*}
which implies the first direction of the statement.

In order to prove the other direction, let us assume that $u$ is harmonic in $\R_+^{n+1}$ and satisfies \refeq{eq:poissonwe}. We first want to show a pointwise estimate for $u$.  So let 
$(x_0,t_0)\in  \R^{n+1}_+$ be arbitrary, $R=t_0/2$,  and $B=B((x_0,t_0), R)$ be the Euclidean ball 
of radius $R$ centered at  $(x_0,t_0)$. Note that for points $(x,t)\in B$ we have 
$t \simeq R$ and  so
$$d\mu_s(x,t) \simeq R^{-(s+1)} dx\,dt$$ on $B$. Hence \eqref{eq:poissonwe} implies that 
$$ |\{(x,t)\in B: |u(x,t)|>\lambda\}|\lesssim \frac{R}{\lambda^p}$$ 
for each $\lambda>0$. 
Since $u$ is harmonic, it follows that
\begin{align*} |u(x_0,t_0)|&\le \frac {1}{|B|}\int_B|u(x,t)|\, dx dt \\ &\le 
\frac {1}{|B|}\bigg( R^{-n/p}|B| + \int_{R^{-n/p}}^\infty |\{(x,t)\in B: |u(x,t)|>\lambda\}|\, d\lambda\biggr)  \\
& \lesssim R^{-n/p}\simeq |t_0|^{-n/p} . 
\end{align*}
This estimate  implies   that for each $\delta>0$ the function $u$ is bounded on $H_\delta:=\{(x,t)\in \R^{n+1}_+: t\ge \delta\}$. If we set $u_\delta=u(\cdot,\delta)$, then the integral representation for bounded harmonic functions \cite[p.~199, Proposition~1]{St}
 shows  that 
\beqla{intrep}
u(x, t+\delta)=(Pu_\delta)(x,t)
\eeq
for $(x,t)\in \R_+^{n+1}$. To obtain a representation for  $u$ itself, we would like to pass to a  limit here along a  sequence $\delta_k\to 0^+$.

For this we fix $r\in (1,p)$, and  write $u=b+v$, where 
$$b=\min\{\max \{u, -1\},1\}$$ and $v=u-b$. Note that $|b|\le 1$.    We want to show that for suitable $\delta>0$ we have good $L^r$-norm bounds for $v(\cdot, \delta)$.

For this we consider  the slab 
$$A_k:=\{ (x,t)\in \R_+^{n+1}\midd   2^{-(k+1)}<t<2^{-k}\}.  $$  
for  $k\in \N$. 
Since 
$$ d\mu_s(x,t)\simeq 2^{k(s+1)}\,  dxdt$$ 
on $A_k$, our assumption  \eqref{eq:poissonwe} implies that
$$ |\{(x,t)\in A_k: |u(x,t)|>\lambda\}|\lesssim  \frac{2^{-k}}{\lambda^p}
$$ for $\lambda>0$.  If $|v|>\lambda>0$, then $|u|> 1+\lambda$, and so
\begin{align*} \int_{A_k} |v(x,t)|^r\, dxdt&=r\int_0^\infty\lambda^{r-1}|\{(x,t)\in A_k: |v(x,t)|>\lambda\}|\, d\lambda\\
&\lesssim 2^{-k}\int_0^\infty \frac {\lambda^{r-1}}{(1+\lambda)^p }\, d\lambda\simeq 2^{-k}. 
\end{align*} Using   Fubini's theorem we can find $\delta_k>0$ with $2^{-(k+1)}<\delta_k<2^{-k}$
such that for $v_k:=v(\cdot, \delta_k)$ we have $\Vert v_k\Vert_{L^r}\lesssim 1$.  Let $b_k=b(\cdot, \delta_k)$. Then $|b_k|\le 1$. So the  sequences $\{v_k\}$  and  $\{b_k\}$ are  uniformly bounded in $L^r(\R^n)$ and $L^\infty(\R^n)$, respectively.  By passing to a subsequence if necessary, we may assume that $v_k\to f_1\in L^r(\R^n) $ and $b_k\to f_2\in L^\infty(\R^n)$ as $k\to \infty$, where convergence is with respect to the respective  weak-$\ast$ topologies. If $t>0$, then for the Poisson kernel 
$P_t$ we have $P_t\in L^q(\R^n)$ for all $q\in [1,\infty]$ \cite[p.~62]{St}. It follows that 
$(P_t\ast v_k)(x)\to (P_t\ast f_1)(x)$ and $(P_t\ast b_k)(x)\to  (P_t\ast f_1)(x)$
as $k\to \infty$ for all $x\in \R^n$ and $t>0$. 
We now define $f=f_1+f_2\in L^r(\R^n)+L^\infty(\R^n)$. 
Note that $u_{\delta_k}=u(\cdot, \delta_k)=v_k+b_k$. Hence  for each $(x,t)\in \R^{n+1}_+$ we have
\begin{align*}
u(x,t)&= \lim_{k\to \infty} u(x,t+\delta_k)=\lim_{k\to \infty}(Pu_{\delta_k})(x,t)\\
&=\lim_{k\to \infty} (P_t\ast(v_k+b_k))(x)= \lim_{k\to \infty} (P_t\ast v_k)(x)+
\lim_{k\to \infty} (P_t\ast b_k)(x)\\ &=  (P_t\ast f_2)(x)+  (P_t\ast f_1)(x)=(Pf)(x,t).
\end{align*}
In other words, $u=Pf$ is the Poisson extension of $f$.

We want to show that  $f\in L^p(\R^n).$ Since $f\in L^r(\R^n)+L^\infty(\R^n)$, where $r>1$, 
 for almost every $x\in \R^n$ we have \cite[p.~62, Theorem~1]{St}
 $$f(x)=\lim_{t\to 0^+} (Pf)(x,t)=\lim_{t\to 0^+} u(x,t). $$ In particular, 
  since $u$ is continuous,  we can find a measurable function $\varepsilon\: \R^n\to [0,\infty]$ such  that $\varepsilon(x) >0$    and 
$$
|u(x,t)|\geq |f(x)|/2\quad {\rm for \;}\; t\in (0,\varepsilon (x))
$$
for almost every $x\in \R^n$. 
Here we may assume that $\eps(x)=\infty$ if $f(x)=0$. 

For $\lambda>0$ consider the set 
$$S_\lambda:=\{ x\in \R^n:  2^{-p} |f(x)|^p\ge 2\lambda^p\eps(x)^{-s}\}$$
(with the understanding that $\eps(x)^{-s}=0$ if $\eps(x)=\infty$). 
Since $\eps(x)>0$ for almost every $x\in \R^n$, the set  $S_\lambda$ increases to a full measure set as $\lambda\to 0^+$. If $x\in S_\lambda$, then  
 $|t^{s/p}u(x,t)|> \lambda$ for  $(2\lambda/|f(x)|)^{p/s}<t<\varepsilon (x)$. We  deduce that
\begin{align*}\lambda^p \mu_s \{(x,t)\in \R^{n+1}_+ &: |t^{s/p}u(x,t)|> \lambda \}   \\ 
&\geq \lambda^p\int_{S_\lambda}
\left(\int_{(2\lambda/|f(x)|)^{p/s}}^{\eps(x)}t^{-(s+1)}dt \right)\,dx\\
&= \frac1s \int_{ S_\lambda}(2^{-p}|f(x)|^p -\lambda^p\eps(x)^{-s})\,dx\\
& \ge \frac{2^{-p}}{2s}  \int_{S_\lambda}|f(x)|^p\, dx. 
\end{align*}
By letting $\lambda\to 0^+$ we conclude 
$$
 \|f\|_{L^p}^p \lesssim \| t^{s/p}u(x,t)\|^p_{L^{p,\infty}(\R^{n+1}_+,\, \mu_s)}< \infty. 
$$ The claim follows. 
\end{proof}

In the next corollary we consider  locally integrable functions $f$ on $\R^n$ that have  distributional partial derivatives $\partial_i f$  in $ L^p(\R^n)$ for $i=1, \dots, n$, where $p>1$. We denote by $\nabla f=(\partial_1f, \dots,\partial_nf)$ the gradient of such a function $f$, and by $|\nabla f|$ the pointwise Euclidean norm of $\nabla f$. 

\begin{corollary}\label{co:sobolevweak}
Let $p>1$, $s>0$, and  $u=Pf$, where $f\in L^p(\R^n)$.
Then   $|\nabla f|\in  L^p (\R^n)$ if and only if 
\beqla{eq:poissonweak1}
t^{s/p}|\nabla u(x,t)|\in L^{p,\infty}\big ( \R_+^{n+1},t^{-(s+1)}dx\,dt).
\eeq
\end{corollary}

\begin{proof}
We consider $u=Pf$, where $f\in L^p(\R^n)$, and use the notation $u_t=u(\cdot,t)$ for $t>0$. 

 We first assume  that $f$ has  distributional partial derivatives $\partial_if\in L^p(\R^n)$ for $i=1, \dots, n$.    Then for each $t>0$ we have 
 $$\partial_i u_t =\partial_i (P_t\ast f)= P_t \ast \partial_if. $$ 
 This shows that all partial derivatives  $\partial_i u$ of $u$ in the $x$-direction are Poisson extensions of functions in $L^p(\R^n)$. To get a similar statement also for the partial derivative $\partial_tu$ in the $t$-direction, we use 
 the  Riesz transforms $R_i$, $i=1, \dots, n$ (see \cite[p.~57]{St} for the definition).  Since  $p>1$, these are bounded operators on $L^p(\R^n)$  \cite[Chapter~3.1]{St},  and so  we have 
$R_i(\partial_if)\in L^p(\R^n)$. Moreover, if we define 
$$g=-\sum_{i=1}^n R_i(\partial_i f)\in  L^p(\R^n), $$
then $\partial_tu_t =P_t\ast g$.  This can  easily verified 
if $f$ is $C^\infty$-smooth and has compact support. 
The general case  follows from the density of such functions in the Sobolev space 
 $W^{1,p}(\R^n)$.
We conclude that all partial derivatives of $u$ are Poisson extensions of functions in $L^p(\R^n)$. 
Condition~\eqref{eq:poissonweak1} now follows from Proposition~\ref{prop:fuweak}.
 
For the converse direction suppose that $u=Pf$ satisfies  \eqref{eq:poissonweak1}.
Then by  Proposition~\ref{prop:fuweak} we know that for  suitable functions $g_1, \dots, 
g_n\in L^p(\R^n)$ we have 
$\partial_iu_t=P_t\ast g_i $ for $i=1, \dots, n$ and all $t>0$. Now standard properties of the Poisson kernel give the convergence  $u_t=P_t\ast f\to f$ and  $\partial_iu_t=P_t\ast g_i \to g_i$ in $L^p(\R^n)$ as $t\to 0^+$ \cite[p.~62, Theorem~1]{St}. This implies that
 $g_i$ is the distributional partial derivative $\partial_if$ of $f$. 
So  $\partial_if=g_i\in L^p(\R^n)$ for $i=1, \dots,n$ and  the claim follows. 
\end{proof}

\begin{remark}{\rm
The previous  corollary essentially gives a characterization of 
functions  $f$  in the Sobolev space $W^{1,p}(\R^n)$ in terms of their Poisson 
extension $u$. For simplicity we made the a priori assumption  $f\in L^p(\R^n)$.  One can relax this integrability condition on $f$ 
and   prove a  more general  result characterizing functions 
 $f$ in the homogeneous Sobolev space $\dot W^{1,p}(\R^n)$ consisting
 of   locally integrable function $f$ on $\R^n$  with distributional derivatives 
  in $ L^p (\R^n)$. 
  For this one first  checks that the Poisson extension is well-defined for 
   each   $f\in \dot W^{1,p}(\R^n)$. Then a generalization of  the above corollary reads as follows:  a locally integrable function $f$ on $\R^n$ belongs to the homogeneous Sobolev space
$\dot W^{1,p}(\R^n)$ if and only if it satisfies the integrability condition $\int_{\R^n}|f(x)|(1+|x|)^{-(n+1)}\, dx<\infty$, and \eqref{eq:poissonweak1} is valid for 
its Poisson extension $u$. }
\end{remark}

The upper half-space $\R^{n+1}_+$ carries the usual   hyperbolic metric. This  is the Riemannian metric  obtained from rescaling  the Euclidean metric  by the factor $1/t$ at the point $(x,t)$. Corresponding to this metric,  one  has  a {\em hyperbolic gradient} $\nabla_hu$ of a smooth function $u$ on $\R^{n+1}_+$ defined by $\nabla_hu(x,t)=t\nabla u(x,t)$ for $(x,t)\in \R^{n+1}_+$, and 
 {\em hyperbolic measure} $\mu_h$ given by   $d\mu_h=t^{-(n+1)}dxdt$.
By specializing   to the case $s=p=n\ge 2$ in Corollary~\ref{co:sobolevweak},  we obtain the following analog of Theorem~\ref{main}.  
\begin{corollary}\label{co:Q-sobolevweak} Let $n\ge 2$ and suppose 
$u=Pf$, where $f\in L^n(\R^n)$. Then $|\nabla f|\in L^n(\R^n)$ if and only if 
$
|\nabla_hu|\in L^{n,\infty}\big ( \R_+^{n+1},d\mu_h). 
$
\end{corollary}

Most of  the results in this section 
remain valid if the Poisson extension  is replaced by other  convolution approximations. We will not pursue this in detail, but 
limit ourselves to recording the following key fact whose proof we will only sketch.

\begin{lemma}\label{le:otherkernel} Let $p>1$, $f\in W^{1,p}(\R^n)$, and 
 $k$ be a bounded and non-negative kernel on $\R^n$
with 
\begin{equation}\label{eq:decay}
0 \; < \; \int_{\R^n}(1+|x|)k(x)\, dx \; <\; \infty.
\end{equation}
Define $\displaystyle u(x,t)=\int_{\R^n}   f(x-y) k(y/t) t^{-n}\, dx$ 
 for $(x,t)\in \R^{n+1}_+$. Then 
 
\begin{equation}\label{eq:extension}
\| f\|_{W^{1,p}}  \simeq \| f\|_{L^p} +\textnormal{ess \!sup}_{t>0}
\Vert |\nabla u|(\cdot, t)\Vert_{L^p}.
%
%
%
\end{equation}
\end{lemma}

Here $\nabla u$ has to be interpreted as the distributional gradient of $u$, and 
we take the essential supremum of the $L^p$-norm of $|\nabla u|(\cdot, t)$ 
for $t>0$. 

Note that condition \eqref{eq:decay} is not true for the Poisson kernel $k=P_1$.

\begin{proof}[Outline of proof] The statement follows from standard approximation properties of convolutions. The only  issue that is not entirely straightforward is how  to control 
the (distributional)  $t$-derivative $\partial_t u$ in terms of the  derivatives $\partial_i f$ of our given function $f\in W^{1,p}(\R^n)$.  

For this we first  assume  in addition to our hypotheses  that $k$ is smooth and compactly supported.   
  Then computation shows that for $t>0$ and $y=(y_1, \dots, y_n)\in \R^n$ we have  
\begin{eqnarray}
\frac{\partial }{\partial t}\big(t^{-n}k(y/t)\big)&=& -nt^{-(n+1)}k(y/t)-t^{-(n+2)}\sum_{i=1}^ny_i(\partial_ik)(y/t)\nonumber\\
&=&\; -t^{-(n+1)}\sum_{i=1}^n\frac{\partial }{\partial y_i}\big(y_i k(y/t)).\nonumber
\end{eqnarray}
Integration by parts then gives    
\begin{equation*}
\partial_tu(x,t)=\sum_{i=1}^n\int_{\R^n} \partial_i f(x-y) v_i(y/t)t^{-n}\, dy,
\end{equation*}
where $v_i(y):=y_ik(y).$  Our hypothesis \eqref{eq:decay} implies  $k\in L^1(\R^n)$ and  $v_i\in L^1(\R^n)$ for $i=1, \dots, n$. 
An $L^p$-bound for $\partial_tu(\cdot,t)$ as desired immediately  follows.

For  general kernels $k$ one uses an approximation argument 
  to reduce  to the case of  smooth and compactly supported kernels.
\end{proof}

By using
the  previous   lemma, one can  (for example)  replace the Poisson extension in Corollary~\ref{co:sobolevweak} by the ball averages 
$$u(x,t)=\frac{1}{|B(x,t)|}\int_{B(x,t)}f(y)\,dy. $$ In the proof one uses 
 standard approximation properties of convolutions combined with 
  the fact that
 the radial maximal function $\sup_{t>0}u(\cdot,t)$ is dominated by  the Hardy-Littlewood maximal function. This last statement is true in greater generality 
  if one demands, in addition to \eqref{eq:decay},
that $k$ is a radially decaying function. For kernels that do not satisfy \eqref{eq:decay} such as the  Poisson kernel,
the issues become more involved.

\end{document}